\documentclass[leqno,12pt]{article} 
\usepackage{url,latexsym,amsmath, amssymb}
\usepackage{amsthm}
\theoremstyle{plain}
\newtheorem{theorem}{\indent\sc Theorem}[section]
\newtheorem{lemma}[theorem]{\indent\sc Lemma}
\newtheorem{corollary}[theorem]{\indent\sc Corollary}
\newtheorem{proposition}[theorem]{\indent\sc Proposition}

\theoremstyle{definition}

\newtheorem*{remark0}{\indent\sc Remark}
\newenvironment{prfof}[1]{\begin{trivlist}\item[\indent \sc{Proof of #1.}]}{
  \end{trivlist} \medskip \par}
\def\prpb{\begin{proposition}}\def\prpe{\end{proposition}}
\def\lemb{\begin{lemma}}\def\leme{\end{lemma}}
\def\thmb{\begin{theorem}}\def\thme{\end{theorem}}
\def\corb{\begin{corollary}}\def\core{\end{corollary}}
\def\prfb{\begin{proof}}\def\prfe{\end{proof}}
\def\prfofb#1{\begin{prfof}{#1}}\def\prfofe{\end{prfof}}
\def\remb{\begin{remark0}}\def\reme{\end{remark0}}
\def\prpa#1{\label{p:#1}}\def\prpu#1{Proposition~{\rm \ref{p:#1}}}
\def\lema#1{\label{l:#1}}\def\lemu#1{Lemma~{\rm \ref{l:#1}}}
\def\thma#1{\label{t:#1}}\def\thmu#1{Theorem~{\rm \ref{t:#1}}}
\def\cora#1{\label{c:#1}}\def\coru#1{Corollary~{\rm \ref{c:#1}}}
\def\seca#1{\label{s:#1}}\def\secu#1{Section~{\rm \ref{s:#1}}}

\def\itmb{\begin{enumerate}}\def\itme{\end{enumerate}}
\def\ittb{\begin{description}}\def\itte{\end{description}}
\def\eqnb{\begin{equation}}
\def\arrb#1{\begin{array}{#1}}\def\arre{\end{array}}
\def\eqnu#1{{\rm (\ref{e:#1})}}

\def\QED{\relax\ifmmode\let\@tempa\relax\ifcase\@eqcnt\def\@tempa{& & &}\or
\def\@tempa{& &}\else\def\@tempa{&}\fi\@tempa $\Box$ \else\hfill $\Box$ \fi}
\def\DDD{\relax\ifmmode\let\@tempa\relax\ifcase\@eqcnt\def\@tempa{& & &}\or
  \def\@tempa{& &}\else\def\@tempa{&}\fi\@tempa $\Diamond$
  \else\hfill $\Diamond$ \fi}
\def\dsp{\displaystyle}
\def\eps{\epsilon}
\def\limf#1{\displaystyle \lim_{#1\to\infty}}

\def\limsup{\displaystyle \mathop{\overline{\lim}}\limits}

\def\diff#1#2{\dsp\frac{d\,#1}{d#2}}
\def\pderiv#1#2{\dsp\frac{\partial\,#1}{\partial#2}}
\def\le{\leqq} \def\ge{\geqq} 
\def\reals{{\mathbb R}}
\def\nreals#1{{\mathbb R}^{#1}}\def\preals{{\mathbb R_+}}

\def\nintegers{{\mathbb N}}\def\rationals{{\mathbb Q}}
\def\prb#1{\def\prbone{#1}
  \ifx\prbone\empty{\mathrm{P}}\else{\mathrm{P[\;}}#1{\mathrm{\;]}}\fi}
\def\prbseq#1#2{\def\prbseqone{#2}
  \ifx\prbseqone\empty{\mathrm{P}}_{#1}\ignorespaces
  \else{\mathrm{P}}_{#1}{\mathrm{[\;}}#2{\mathrm{\;]}}\fi}
\def\EE#1{{\mathrm{E[\;}}#1{\mathrm{\;]}}}
\def\EEseq#1#2{\def\EEseqone{#2}
  \ifx\EEseqone\empty{\mathrm{E}}_{#1}\else
 {\mathrm{E}}_{#1}{\dsp\mathrm{[\;}}#2{\mathrm{\;]}}\fi}
\def\VV#1{{\mathrm{V[\;}}#1{\mathrm{\;]}}}
\def\VVseq#1#2{\def\VVseqone{#2}
  \ifx\VVseqone\empty{\matrm{V}}_{#1}\else
 {\mathrm{V}}_{#1}{\dsp\mathrm{[\;}}#2{\mathrm{\;]}}\fi}

\def\ssN{^{(N)}}
\def\ssNn{^{(N,n)}}

\def\chrfcn#1{\mathop{\mathbf{1}}\nolimits_{#1}}
\def\suppt{\mathop{\mathrm{suppt}}\nolimits}
\def\msp{{\cal M}(\preals)}
\def\fsp{L^1_{{\rm loc}}(\preals)}

\makeatletter
\def\address#1#2{\begingroup
\noindent\parbox[t]{7.8cm}{%
\small{\scshape\ignorespaces#1}\par\vskip1ex
\noindent\small{\itshape E-mail address}%
\/: #2\par\vskip4ex}\hfill%
\endgroup}%
\makeatother
\title{\uppercase{
 Stochastic ranking process with time dependent intensities
}}
\author{
\textsc{
Yuu Hariya,
Kumiko Hattori$^{*}$,
} \and \textsc{
Tetsuya Hattori$^{**}$,
Yukio Nagahata,
} \and \textsc{
Yuusuke Takeshima and
Takahisa Kobayashi
}
}
\date{}
\begin{document}
\maketitle
\thispagestyle{empty}

\footnote{ 
$^{*}$
Partly supported by the Grant-in-Aid for Scientific Research (C),
Japan Society for the Promotion of Science. 

$^{**}$
Partly supported by the Grant-in-Aid for Scientific Research (B),
Japan Society for the Promotion of Science. 

}

\begin{abstract}
We consider the stochastic ranking process with 
the jump times of the particles determined by
Poisson random measures.
We prove that the joint empirical distribution of scaled position and
intensity measure converges almost surely in the infinite particle limit.
We give an explicit formula for the limit distribution and show that
the limit distribution function is a unique global classical solution
to an initial value problem for a system of a first order non-linear 
partial differential equations with time dependent coefficients.
\end{abstract}

\section{Introduction.}
\seca{1}
Let $\msp$ be the space of Radon measures $\rho$ on the Borel $\sigma$-algebra
${\cal B}(\preals)$ of non-negative reals $\preals$.
Let $N$ be a positive integer, and let $\nu\ssN_i$, $i=1,2,\ldots,N$, 
be independent Poisson random measures (Poisson point processes) 
on $\preals$\,,
defined on a probability space $(\prb{},{\cal F},\Omega)$.
For each $i$, denote the intensity measure of $\nu\ssN_i$ by $\rho\ssN_i$;
\eqnb
\label{e:PRM}
\EE{\nu\ssN_i(A)}=\rho\ssN_i(A),\ A\in{\cal B}(\preals).
\end{equation}
Throughout this paper,
we assume $\rho\ssN_i\in\msp$ and that $\rho\ssN_i$ is continuous
(i.e., $\rho\ssN_i(\{t\})=0$, $t\ge0$) for all $N$ and $i$.

Let $x\ssN_1,x\ssN_2,\ldots,x\ssN_N$ be a permutation of $1,2,\ldots,N$,
and define a process 
\[X\ssN=(X\ssN_1,\ldots,X\ssN_N)\]
 by
\eqnb
\label{e:SIformSRP}
\arrb{l} \dsp
X\ssN_i (t) 
= x\ssN_i + \sum_{k=1}^N \int_0^t
 \chrfcn{X\ssN_k(s-0)>X\ssN_i(s-0)}\, \nu\ssN_k(ds)
+ \int_0^t (1-X\ssN_i(s-0))\, \nu\ssN_i(ds),
\\ \ i=1,2,\ldots,N,\ t\ge 0\,,
\arre
\end{equation}
where, $\chrfcn{A}$ is the indicator function of event $A$.

Denote the unit measure concentrated on $c$ by $\delta_c$\,.
With probability $1$ we can write
\eqnb
\label{e:PoissonUnitMeasRep}
\nu\ssN_i=\sum_{j=1}^{\infty} \delta_{\tau\ssN_{i,j}},
\ i=1,2,\ldots,N,
\end{equation}
where, with probability $1$,
$\tau\ssN_{i,j}$'s are random variables satisfying
$0<\tau\ssN_{i,1}<\tau\ssN_{i,2}<\cdots$, $i=1,2,\ldots,N$, and
$\tau\ssN_{i,j}\ne \tau\ssN_{i',j'}$ if $(i,j)\ne (i',j')$.
In the following, we work on the event that these inequalities hold.

The right hand side of \eqnu{SIformSRP} is a simple function in $t$.
At $t=\tau_{i,j}$ we see
\[
 X\ssN_i (\tau_{i,j})-X\ssN_i(\tau_{i,j}-)= 1-X\ssN_i(\tau_{i,j}-) ,
\]
which implies 
\eqnb
\label{e:defprocjx}
X\ssN_i (\tau_{i,j}) =1.
\end{equation}
With similar consideration, we see that the process $X\ssN$ is uniquely
determined by \eqnu{SIformSRP}: Explicitly, we have, for $i=1,\ldots,N$,
\eqnb
\label{e:defprocjx1}
X\ssN_i(t)=\left\{\arrb{ll} \dsp
x\ssN_{i}+\sum_{i';\; x\ssN_{i'}>x\ssN_{i}}
 \chrfcn{\tau\ssN_{i',1} \le t}
&\ \ \  0\le t< \tau\ssN_{i,1}\,, 
\\ \dsp
1+\sum_{i'=1}^N
 \chrfcn{\exists j'\in\nintegers;\; \tau\ssN_{i,j}<\tau\ssN_{i',j'}\le t}
&\ \ \  \tau\ssN_{i,j}\le t<\tau\ssN_{i,j+1},\ j=1,2,3,\ldots.
\arre\right.
\end{equation}

In the case of the (homogeneous) Poisson process 
(i.e., the case $\rho\ssN_i((0,t])=w\ssN_i t$,
 $t\ge0$, for positive constants $w\ssN_i$),
a discrete time version of the process \eqnu{defprocjx1}
has been known for a long time 
\cite{Tsetlin1963,mv2frnt1,mv2frnt2,mv2frnt3,Letac1974,Kingman75} 
and is called move-to-front (MTF) rules.
The process has, in particular, been extensively studied as a model of 
least-recently-used (LRU) caching in the field of information theory 
\cite{Rivest1976,Fagin77,Bitner79,CHS88,BlomHolst91,Rodrigues95a,Fill96JTP,%
mv2frnt4,Fill96TCS,Jelenkovic99,Jelenkovic03},
and also is noted as a time-reversed process of top-to-random shuffling.
With a great advance in the internet technologies,
a new application of the process appeared \cite{HH072,mv2frnt}.
The ranking numbers such as those found in the web pages of
online bookstores are found to follow the predictions of the model.

In \cite{HH071}, the case where $\nu\ssN_i$'s are
(homogeneous) Poisson processes with $\rho\ssN_i((0,t])=w\ssN_i t$
is considered, and the joint empirical distribution of jump rate $w\ssN_i$ 
and normalized position
\eqnb
\label{e:YN}
 Y\ssN_i(t)=\frac1N(X\ssN_i(t)-1),
\end{equation}
given by
$\dsp
\mu\ssN_{t} =\frac1N \sum_{i=1}^N \delta_{(w\ssN_{i},Y\ssN_i(t))}\,,
$
is studied. (We will abuse notation slightly and denote a unit measure
on any space by $\delta_c$.)
It is proved in \cite{HH071} that a scaling limit
\eqnb
\label{e:mainlimit}
\mu_{t}= \limf{N}\mu\ssN_{t}
\end{equation}
exists (under reasonable assumptions), 
and an explicit formula for $\mu_t$, which is a 
deterministic distribution on $\preals\times [0,1)$, is given.
In \cite{HH072}, it is proved that,
if the scaling limit of the jump rate distribution is a discrete distribution,
the limit $\mu_t$ is the unique time global solution to 
an initial value problem for a system of 
first order non-linear partial differential equations
(inviscid Burgers equations with a term representing evaporation).
The structure of the explicit formula for $\mu_t$ is
naturally explained by a standard method of characteristic curves
for the solution to the partial differential equations.

In the present paper, 
we will generalize the main results of \cite{HH071,HH072} to the case
where $\nu\ssN_i$'s are
Poisson random measures.
We shall call the process $X\ssN$ defined by \eqnu{SIformSRP}, 
or equivalently by \eqnu{defprocjx1},
a stochastic ranking processes after \cite{HH071,HH072,HH073}.

Put
\eqnb
\label{e:XCN}
X\ssN_C(t)= \sum_{i=1}^N \chrfcn{\tau\ssN_{i,1}\le t},\ \ t\ge0.
\end{equation}
$X\ssN_C(t)$ is a random variable which denotes the position of the
boundary between the top side $x\le X\ssN_C(t)$
and the tail side $x> X\ssN_C(t)$,
where each particle in the top side 
(i.e., $i$ which satisfies $X\ssN_i(t)\le X\ssN_C(t)$)
has experienced jump to the top by time $t$ (i.e., $\tau\ssN_{i,1}\le t$),
and the particles in the tail side are those particles which have not
jumped to the top by time $t$.
\prpb
\prpa{yC}
Let $t\ge0$, and assume that a sequence of distributions 
$\{\lambda\ssN_t\;;\ N\in\nintegers\}$ on $\preals$ defined by
\eqnb
\label{e:lambdaN}
\lambda\ssN_t=\frac1N \sum_{i=1}^N \delta_{\rho\ssN_i((0,t])}
\end{equation}
converges weakly as $N\to\infty$ to a probability distribution $\lambda_t$\,.
Then the scaled position of the boundary
\eqnb
\label{e:YCN}
Y\ssN_C(t)=\frac1N X\ssN_C(t)=\frac1N \sum_{i=1}^N\chrfcn{\tau\ssN_{i,1}\le t}
\end{equation}
converges almost surely as $N\to\infty$ to
\eqnb
\label{e:yC}
y_C(t)=1-\int_0^{\infty} e^{-s} \lambda_t(ds).
\end{equation}
\DDD\prpe
\prfb
The definition \eqnu{YCN} implies that 
$Y\ssN_C(t)-\EE{Y\ssN_C(t)}$
is an arithmetic mean of independent variables 
\[ Z\ssN_i=\chrfcn{\tau\ssN_{i,1}\le t}-\prb{\tau\ssN_{i,1}\le t},
\ \ i=1,2,\ldots,N, \]
with bounded $4$th order moment.
(In fact, $|Z\ssN_i|\le 1$, for all $N$ and $i$.)
Hence, 
\[
 \EE{\sum_{N=1}^{\infty} (Y\ssN_C(t)-\EE{Y\ssN_C(t)})^4}
= \sum_{N=1}^{\infty} \EE{(Y\ssN_C(t)-\EE{Y\ssN_C(t)})^4}<\infty, \]
which implies 
\[
Y\ssN_C(t)-\EE{Y\ssN_C(t)}\to 0,\ \mbox{ a.e., as }\ N\to\infty.
\]
On the other hand, definition of Poisson random measure implies
\[
\EE{Y\ssN_C(t)}=\frac1N\sum_{i=1}^N \prb{\tau\ssN_{i,1}\le t}
= \frac1N\sum_{i=1}^N (1-e^{-\rho\ssN_i((0,t])})
= 1- \int_0^{\infty} e^{-s} \lambda\ssN_t(ds),
\]
which converges to \eqnu{yC} by assumption.
\prfe

Since by \prpu{yC} we have almost sure convergence at each time $t$, 
we have almost sure convergence for all rational number times simultaneously.
By definition, $y_C(t)$ and $Y\ssN_C(\omega)(t)$, $\omega\in\Omega$, are 
non-decreasing in $t$. Hence, if $y_C(t)$ is continuous, we have almost sure 
convergence as a function in $t$.
\corb
\cora{yC}
In addition to the assumptions in \prpu{yC}, assume that
$\lambda_t$ is continuous in $t$ with respect to the 
topology of weak convergence.  
Then for almost all sample $\omega\in\Omega$,
$Y\ssN_C(\omega):\ \preals\to[0,1)$ defined by \eqnu{YCN}
converges pointwise in $t$ as $N\to\infty$ to a deterministic function
$y_C:\ \preals\to[0,1)$ defined by \eqnu{yC}.
\DDD\core

\prpu{yC} is a generalization to inhomogeneous case 
of \cite[Proposition 2]{HH071} for the (homogeneous) Poisson process.
The correspondence with $\lambda_t$ in \prpu{yC} and $\lambda$
in \cite{HH071} is given by $\lambda_t((0,c\,t])=\lambda((0,c])$.
\eqnu{lambdaN} implies that $\lambda_t$ is the asymptotic distribution
of the expectation of number of jumps to rank $1$ for each particle in the
time interval $(0,t]$.

Consider a joint empirical distribution $\mu\ssN$ of intensity measure
$\rho\ssN_i$ and scaled position $Y\ssN_i$ of the stochastic ranking 
process:
\eqnb
\label{e:muN}
\mu\ssN_{t} =\frac1N \sum_{i=1}^N \delta_{(\rho\ssN_i,Y\ssN_i(t))}\,,
\ \ t\ge0.
\end{equation}
$\mu\ssN_t$\,, $N\in\nintegers$, are random variables whose samples
are distributions on
the product space $\msp\times [0,1)$ of space of Radon 
measures $\msp$ and an interval $[0,1)\subset \preals$.

We consider the standard vague topology on $\msp$, 
that is, a sequence $\{\rho_n\}\subset \msp$ converges to $\rho\in\msp$
if and only if
\eqnb
\label{e:vagueconv} \ \ \ 
\limf{n}\int_{\preals}f(s)\,\rho_n(ds)=\int_{\preals}f(s)\,\rho(ds),
\end{equation}
for all continuous function $f$ with compact support.
Since $\preals$ is a Polish space, i.e., 
complete and separable metric space,
so is $\msp$ \cite[Theorem 31.5]{Bauer}, and consequently,
$\msp\times[0,1)$ is also a 
Polish space \cite[Example 26.2]{Bauer}. 

Assume that a sequence of initial configurations
\[
\mu\ssN_{0} =\frac1N \sum_{i=1}^N \delta_{(\rho\ssN_i,(x\ssN_i-1)/N)}\,,
\ \ N=1,2,\ldots,
\]
converges weakly as $N\to\infty$ to a probability distribution $\mu_0$
on $\msp\times [0,1)$. 
Then, in particular,
\eqnb
\label{e:LambdaN}
\arrb{l}\dsp
\Lambda\ssN(d\rho):= \mu\ssN_{0}(d\rho\times[0,1)) 
=\frac1N \sum_{i=1}^N \delta_{\rho\ssN_i}(d\rho)
\ \to\ \Lambda(d\rho):=\mu_0(d\rho\times[0,1)),
\\
\mbox{ weakly, as } N\to\infty.
\arre
\end{equation}
Note also that $\lambda\ssN_t$ in \eqnu{lambdaN} has an expression
\eqnb
\label{e:lambdaNLambda}
\lambda\ssN_t=\int_{\msp} \delta_{\rho((0,t])} \Lambda\ssN(d\rho).
\end{equation}
We shall generalize \eqnu{lambdaNLambda} and define, for $0\le s\le t$,
\eqnb
\label{e:lambdaNstLambda}
\lambda\ssN_{s,t}=\int_{\msp} \delta_{\rho((s,t])} \Lambda\ssN(d\rho).
\end{equation}
\thmb
\thma{HDLinhomog}
Assume that $\mu\ssN_0\to \mu_0$ weakly as $N\to\infty$
for a probability distribution $\mu_0$ on $\msp\times [0,1)$. 
Assume that for each $(s,t)$ satisfying $t\ge s\ge0$, 
\eqnb
\label{e:lambdaN2lambda}
\lambda\ssN_{s,t}\ \to\ 
\lambda_{s,t}:=\int_{\msp} \delta_{\rho((s,t])} \Lambda(d\rho),
\ \mbox{ weakly as }\ N\to\infty,
\end{equation}
where $\Lambda$ is as in \eqnu{LambdaN}.
Then for any $t>0$, and for almost all sample $\omega\in\Omega$,
the distribution $\mu\ssN_{t}(\omega)$ converges weakly
to a non-random probability distribution $\mu_t$ on $\msp\times [0,1)$.

$\mu_t$ has a following expression in terms of 
$U(d\rho,y,t):=\mu_{t}(d\rho \times [y,1))$.
\eqnb
\label{e:HDLdistri}
U(d\rho,y,t):=\mu_{t}(d\rho \times [y,1))=
\left\{ \arrb{ll}\dsp
e^{-\rho((t-t_0(y,t),t])}\, \Lambda(d\rho)\ \ \  & 0\le y\le y_C(t), \\ \dsp
e^{-\rho((0,t])}\, U(d\rho,\hat{y}(y,t),0)\ \ \  & y_C(t)\le y<1.
\arre \right. 
\end{equation}
Here,
$t_0(y,t)$ is the inverse function with respect to $t_0$ of
\eqnb
\label{e:yA}
 y_A(t_0,t)=1-\int_{\msp} e^{-\rho((t-t_0,t])}\, \Lambda(d\rho),
\ \ 0\le t_0\le t,
\end{equation}
namely,
\eqnb
\label{e:t0}
 t_0(y,t)=\inf\{s\in[0,t] \;;\ y_A(s,t)\ge y \},
\end{equation}
and $\hat{y}(y,t)$ is the inverse function with respect to $y$ of 
\eqnb
\label{e:yB}
 y_B(y,t)=1-\int_{\msp} e^{-\rho((0,t])}\, \mu_{0}(d\rho\times [y,1)),
\ \ t\ge 0,\ 0\le y<1,
\end{equation}
namely,
\eqnb
\label{e:haty}
 \hat{y}(y,t)=\inf\{x\in[0,1) \;;\ y_B(x,t)\ge y \}.
\end{equation}
\DDD\thme
Note that $y_C(t)=y_A(t,t)=y_B(0,t)$.
Note also that, as will be evident from the proof of \thmu{HDLinhomog}
in \secu{HDL} for $0\le y\le y_C(t)$,
the assumption $\mu\ssN_0\to \mu_0$ can be replaced by 
a weaker assumption $\Lambda\ssN_0\to\Lambda$
for $0\le y\le y_C(t)$.

In contrast to \prpu{yC}, we do not have a result analogous to \coru{yC}
for \thmu{HDLinhomog}, because we can expect no monotonicity for $\mu\ssN_t$.
If we impose additional conditions,
we may go further and prove almost sure convergence as sequences of 
processes on a finite time interval $[0,T]$, both for $Y\ssN_C\to y_C$
and $\mu\ssN\to \mu$. See \secu{fnctHDL}
for statement (\thmu{fnctHDLinhomog}) and proof.

The structure of the explicit limit formula \eqnu{HDLdistri},
in particular, the appearance of the inverse functions $t_0$ of $y_A$ and 
$\hat{y}$ of $y_B$\,, can mathematically be understood
through a system of partial differential equations,
which is a generalization of that in \cite{HH072}.
To avoid notational complication, 
consider the case that the limit distribution
$\Lambda$ is supported on a discrete set:
$
 \Lambda= \sum_{\alpha} r_{\alpha} \delta_{\rho_{\alpha}}\,.
$
Then \eqnu{HDLdistri} implies, for
$U_{\alpha}(y,t):= \mu_{t}(\{\rho_{\alpha}\}\times [y,1))$,
\eqnb
\label{e:HDLdistridiscrete}
U_{\alpha}(y,t)=\left\{ \arrb{ll}\dsp
r_{\alpha} \,e^{-\rho_{\alpha}((t-t_0(y,t),t])}\ \ \ & 0\le y\le y_C(t), 
\\ \dsp
U_{\alpha}(\hat{y}(y,t),0)\, e^{-\rho((0,t])}\ \ \ & y_C(t)\le y<1,
\arre \right. 
\end{equation}
where $t_0$ and $\hat{y}$ are inverse functions, respectively, of 
\eqnb
\label{e:BurgersyA}
 y_A(t_0,t)=1-\sum_{\alpha} r_{\alpha} e^{-\rho_{\alpha}((t-t_0,t])},
\end{equation}
and
\eqnb
\label{e:BurgersyB}
 y_B(y,t)=
1-\sum_{\alpha}U_{\alpha}(y,0) e^{-\rho_{\alpha}((0,t])},
\end{equation}
defined by \eqnu{t0} and \eqnu{haty}.
\thmb
\thma{Burgers}
Let $k$ be a positive integer, and for each $\alpha=1,2,\ldots,k$, 
let $r_{\alpha}$ be a positive constant,
$w_{\alpha}:\ \preals\to\preals$ a measurable function
satisfying $w_{\alpha}(t)>0$, $t\ge0$,
and $u_{\alpha}:\ [0,1)\to \preals$
a non-negative smooth strictly decreasing function, satisfying
\eqnb
\label{e:Burgersassump}
 \sum_{\beta=1}^k r_{\beta}=1,\ \ 
 \sum_{\beta=1}^k r_{\beta} w_{\beta}(t) <\infty,\ t\ge0,
\ \mbox{ and }\ 
 \sum_{\beta=1}^k u_{\beta}(y) = 1-y,\ 0\le y<1.
\end{equation}
Then an initial value problem for a system of partial differential equations
\eqnb
\label{e:Burgers}\ \ \ \ \ 
\arrb{l}\dsp
\pderiv{U_{\alpha}}{t}(y,t) 
+\sum_{\beta=1}^k w_{\beta}(t)\, U_{\beta}(y,t)\, \pderiv{U_{\alpha}}{y}(y,t)
=-w_{\alpha}(t) U_{\alpha}(y,t),
\\ \dsp
\  (y,t)\in [0,1)\times \preals,\ \alpha=1,2,\ldots,k,
\arre
\end{equation}
with a boundary condition
\eqnb
\label{e:Burgersboundary}
U_{\alpha}(0,t)=r_{\alpha},\ t\ge0,
\ \alpha=1,2,\ldots,k,
\end{equation}
and initial data 
\eqnb
\label{e:Burgersinitial}
 U_{\alpha}(\cdot,0)=u_{\alpha},
\ \alpha=1,2,\ldots,k,
\end{equation}
has a unique time global classical solution, whose formula is given by
\eqnu{HDLdistridiscrete}
with 
\eqnb
\label{e:Burgersrho}
\rho_{\alpha}((s,t])=\int_s^t w_{\alpha}(u)\,du\ 
\mbox{ and }\ U_{\alpha}(y,0)=u_{\alpha}(y).
\end{equation}
\DDD\thme
As in \cite[\S2]{HH072},
\eqnu{Burgers} is solved by a method of characteristic curves,
and $y_A$, $y_B$, and $y_C$ turn out to be 
the characteristic curves for \eqnu{Burgers},
which mathematically explains how the inverse functions of these functions
appear in the solutions.

For the homogeneous case ($\rho\ssN_i((0,t])=w\ssN_i t$),
\thmu{HDLinhomog} reduces to \cite[Theorem 5]{HH071}
(with slightly weaker assumption on $\mu_0$, $\Lambda$, and $\lambda_t$,
and with stronger convergence in $(\Omega,{\cal F},\prb{})$,
thanks to technical refinement in the proof),
and \thmu{Burgers} reduces to \cite[Theorem 1]{HH072}.
Motivation for extending the previous results to the present case
arises both from mathematical and application point of view.
\ittb
\item[Mathematical: ]
The model is a natural extension of \cite{HH071}, with 
(homogeneous) Poisson processes in the formulation of \cite{HH071}
generalized to (inhomogeneous) Poisson random measures in 
\eqnu{SIformSRP} or \eqnu{defprocjx1}.
Also, as seen from \thmu{Burgers},
the system of PDE corresponding to the limit distribution
is a natural extension of that considered in \cite{HH072}, with
constant coefficients $w_{\alpha}$ in \cite{HH072} generalized to 
time dependent coefficients $w_{\alpha}(t)$ in \eqnu{Burgers}.
On the other hand, 
the space on which $\mu_t$ is defined becomes large;
$\mu_t$ considered in \cite{HH071} is a distribution on $\preals\times [0,1)$,
whereas $\mu_t$ in \thmu{HDLinhomog} is on $\msp\times [0,1)$.
Hence it is necessary to extend the definition of the model,
compared to \cite{HH071,HH072}.

\item[Application: ]
The model has successfully been applied to
statistical explanation of ranking  data at an online bookstore Amazon.co.jp
\cite{HH073,HH072} and data of list of subject titles at a collected 
bulletin board 2ch.net \cite{HH072}.
These data arise as results of social activities,
hence it is inevitable that the data have day-night difference
in their time dependence.
This motivates considering the inhomogeneous cases from an application side.
\itte

Note that we directly see from \eqnu{SIformSRP}, the Markov property
\[ \arrb{l}\dsp
X\ssN_i (t+u) = X\ssN_i(u)
 + \sum_{k=1}^N \int_0^t
 \chrfcn{X\ssN_k(s+u-0)>X\ssN_i(s+u-0)}\, \tilde{\nu}\ssN_k(ds)
\\ \dsp \phantom{X\ssN_i (t+u) = }
+ \int_0^t (1-X\ssN_i(s+u-0))\, \tilde{\nu}\ssN_i(ds), 
\arre \]
where we put $\tilde{\nu}\ssN_i(A)=\nu\ssN_i(A+u)$.
In practical application, this property enables us to shift the time 
origin $t=0$ to the time that a particle we observe jumps to the top,
namely, we may set $X\ssN_i(0)=x\ssN_i=1$, by adjusting the `clock' for
the intensity measure accordingly.
This motivates our formulating the model in terms of Poisson random measures,
even though in \prpu{yC} we apparently do not use Markov properties.

Note also that if $x\ssN_i=1$, then up to the first jump of $i$ to the top,
namely, for $t<\tau\ssN_{i,1}$, comparison of 
\eqnu{defprocjx1} and \eqnu{XCN} leads to
\[ X\ssN_i(t)=X\ssN_C(t)+1, \]
because, if $x\ssN_i=1$ then, $x\ssN_{i'}>x\ssN_i$ for all $i'\ne i$.
Therefore, in practical application, 
we may proceed with observing a trajectory (time development)
of a single particle, putting the time of its first jump to top as $t=0$
and observing until its next jump to top,
and then apply \prpu{yC} or \coru{yC} \cite{HH072,HH073}.

The plan of the paper and a brief description of the role
of the authors are as follows.
In \secu{HDL} we prove \thmu{HDLinhomog},
and we prove \thmu{Burgers} in \secu{Burgers}.
In \secu{fnctHDL}, we state and prove \thmu{fnctHDLinhomog},
time-uniform results corresponding to \prpu{yC} and \thmu{HDLinhomog}.
The core structure of the present work, including basic properties 
of the stochastic ranking process which are essential for the proofs 
of these results, are based on collaboration of K.~Hattori and T.~Hattori.
In extending the previous results 
for the convergence of empirical distribution
on $\preals\times[0,1)$ to $\msp\times [0,1)$, where $\msp$ is a space of 
Borel measures, we have to reformulate the process using
Poisson random measures and provide abstract measure theory result
\lemu{Bauer}, for which collaboration with Hariya is crucial.
Convergence result as measure valued processes developed in \secu{fnctHDL}
is achieved by collaboration with Nagahata.
Also, various technical refinements,
implying in particular  stronger convergence with
less assumptions for the uniform intensity case \cite{HH071}, are results
of the collaboration of these $4$ authors.
In \secu{Appl} we consider a simple case where
the intensities of the Poisson random measures have a common time dependence,
and prove another scaling limit for the particle trajectory,
corresponding to a time change with respect to the intensity.
This is a result of collaboration of T.~Hattori, Hariya, Kobayashi, 
and Takeshima at Tohoku University, and provides
a mathematical result of scaling limit with time changes, as well as
a practically useful formula in applying the present results to 
online rankings.
A practical method based on this mathematical result is partly checked
by actual data obtained at 2ch.net in the master theses of
Kobayashi and Takeshima (unpublished).
In Appendix, we give remarks to be kept in mind when applying our results
to practical data through statistical analysis.

\smallskip\par\textbf{Acknowledgment.}

The authors would like to thank Professor Masayoshi Takeda 
for collaboration at Tohoku University.

\section{Proof of \protect\thmu{HDLinhomog}.}
\seca{HDL}

Throughout this section, 
we assume that the assumptions of \thmu{HDLinhomog} hold.

We first note the following rather technical generality.
\lemb
\lema{Bauer}
Let $t>0$.
If, for each $y\in[0,1)$ and for each
bounded continuous function $g:\ \msp\to \reals$,
there exists $\tilde{\Omega}$ with $\prb{\tilde{\Omega}}=1$ such that
\eqnb
\label{e:HDLproof1}
\arrb{l}\dsp
\limf{N}
\frac1N\sum_{i=1}^N g(\rho\ssN_i)\, \chrfcn{Y\ssN_i(t)\ge y}(\omega)
\\ \dsp {}
= \left\{ \arrb{ll}\dsp
\int_{\msp} g(\rho)\,e^{-\rho((t-t_0(y,t),t])}\, \Lambda(d\rho)
\ \ \ & 0\le y\le y_C(t), \\ \dsp
\int_{\msp} g(\rho)\,e^{-\rho((0,t])}\,\mu_{0}(d\rho\times [\hat{y}(y,t),1))
\ \ \ & y_C(t)\le y<1,
\arre \right. 
\arre
\end{equation}
holds on $\tilde{\Omega}$,
then the claim of \thmu{HDLinhomog} holds for this $t$.
\DDD\leme
The point here is that $\tilde{\Omega}$ may depend on $y$ and $g$,
while \thmu{HDLinhomog} claims the existence of a sample set,
independently of test functions.

We make use of the results in \cite[Exercises 30.3, 31.2]{Bauer}
for a proof of \lemu{Bauer}.
Note that $\msp$ is not locally compact, while local compactness is assumed
in the relevant results of the reference.
We prepare the next Lemma to fill the gap.
\lemb
\lema{HDLproof0}
There exists a \textit{countable} subset 
${\cal T}=\{f_n\;;\ n\in\nintegers\}$ of
uniformly 
continuous functions $f_n:\ \msp\times[0,1)\to\reals$,
such that if for each $f_n\in {\cal T}$
\eqnb
\label{e:HDLproof0}
\limf{N} 
\int_{\msp\times[0,1)} f_n(\rho,y)\,\nu_N(d\rho\times dy)=
\int_{\msp\times[0,1)} f_n(\rho,y)\,\nu(d\rho\times dy)
\end{equation}
holds for a sequence of Borel probability measures $\nu_N$
and a Borel probability measure $\nu$ on $\msp\times[0,1)$, 
then $\nu_N\to \nu$, weakly as $N\to\infty$.
\DDD\leme
\prfb
We noted below \eqnu{vagueconv} that $\msp\times[0,1)$ is a Polish space.
Note also that there exists a coutable set of continuous functions 
$\{e_n:\ \preals\to\reals\;;\ n\in\nintegers\}$ of compact support,
such that 
\eqnb
\label{e:rev0724-1}
\arrb{l}\dsp
d((\rho_1,y_1),(\rho_2,y_2))
=|y_1-y_2|+\sum_{n\in\nintegers} 2^{-n} (1\wedge
 \left|\int_{\preals}e_n(s)\,\rho_1(ds)
-\int_{\preals}e_n(s)\,\rho_2(ds)\right|),
\\ \dsp
\ (\rho_i,y_i)\in\msp\times[0,1),\ i=1,2,
\arre\end{equation}
defines a metric $d$ compatible with the topology we are considering
\cite[(31.4)]{Bauer}.

Denote a set of sequences by $\nreals{\infty}=\{x=(x_1,x_2,\ldots)\}$,
and define a metric $d'$ on $\nreals{\infty}\times[0,1)$ by
\eqnb
\label{e:rev0724-2}
d'((x_1,y_1),(x_2,y_2))=|y_1-y_2|+
\sum_{n\in\nintegers} 2^{-n} (1\wedge |x_{1,n}-x_{2,n}|)
\end{equation}
where $x_i=(x_{i,1},x_{i,2},\ldots)$, $i=1,2$.
We have a natural one-to-one map 
$\iota=(\iota_1,\iota_2,\ldots,\iota_0)
:\ \msp\times[0,1)\to \nreals{\infty}\times[0,1)$
defined by
\eqnb
\label{e:rev0724-3}
\iota(\rho,y)_n=\int_{\preals}e_n(s)\,\rho(ds),\ n\in\nintegers,
\ \mbox{ and }\ \iota(\rho,y)_0=y.
\end{equation}
Put 
\eqnb
\label{e:rev0724-4}
E'=\iota(\msp\times[0,1))\subset\nreals{\infty}\times[0,1).
\end{equation}
Then 
\eqnu{rev0724-1}, \eqnu{rev0724-2} and \eqnu{rev0724-3} imply that
$\iota:\ \msp\times[0,1)\to E'$ is a one-to-one onto isometric map.
Since $\msp\times[0,1)$ is complete, 
$E'$ is a closed set in $\nreals{\infty}\times[0,1)$.

Let $F\subset \msp\times[0,1)$ be a closed set.
Since $\iota$ is isometric, $\iota(F)$ is a closed subset of $E'$,
and since $E'$ is a closed set in $\nreals{\infty}\times[0,1)$,
$\iota(F)$ is a closed set in $\nreals{\infty}\times[0,1)$.
Hence, if a sequence of probability measures $\nu_N\circ \iota^{-1}$
on $\nreals{\infty}\times[0,1)$ converges weakly as $N\to\infty$ to 
$\nu\circ \iota^{-1}$, then
\[ \limsup_{N\to\infty}\nu_N(F)=
\limsup_{N\to\infty}\nu_N\circ\iota^{-1}(\iota(F))\le
\nu\circ\iota^{-1}(\iota(F))=\nu(F), \]
which implies $\nu_N\to\nu$, weakly as $N\to\infty$.
Thus the conclusion of \lemu{HDLproof0} is reduced to
a weak convergence $\nu_N\circ \iota^{-1}\to \nu\circ \iota^{-1}$ on
$\nreals{\infty}\times[0,1)$.

For each $k\in\nintegers$ define a projection to finite dimensional space
$\pi_k:\ \nreals{\infty}\times[0,1)\to\nreals{k}\times[0,1)$ by
\eqnb
\label{e:rev0724-6}
\pi_k(x)=(x_1,x_2,\ldots,x_k,y),
\ x=(x_1,x_2,\ldots,y)\in\nreals{\infty}\times[0,1).
\end{equation}
Then $\nu_N\circ\iota^{-1}\circ\pi_k^{-1}$ and 
$\nu\circ\iota^{-1}\circ\pi_k^{-1}$ are probability measures 
on $\nreals{k}\times[0,1)$.
Note that a Borel probability measure on Polish space is a Radon measure
\cite[Theorem 26.3]{Bauer}, and that the vague convergence of
probability measures to a probability measure on $\nreals{k}$
is equivalent to the weak convergence \cite[Theorem 30.8]{Bauer}.
Since $\nreals{k}\times[0,1)$ is a locally compact Polish space,
there exists a \textit{countable} subset 
${\cal T}_k=\{f_{k,i}\;;\ i\in\nintegers\}$ of
continuous functions $f_{k,i}:\ \nreals{k}\times[0,1)\to\reals$ 
with compact support,
such that if for each $f_{k,i}\in {\cal T}_k$
\eqnb
\label{e:rev0724-7}
\limf{N} 
\int_{\nreals{k}\times[0,1)} 
f_{k,i}(z)\,\nu_N\circ\iota^{-1}\circ\pi_k^{-1}(dz)=
\int_{\nreals{k}\times[0,1)} f_{k,i}(z)\,\nu\circ\iota^{-1}\circ\pi_k^{-1}(dz)
\end{equation}
holds, then $\nu_N\circ\iota^{-1}\circ\pi_k^{-1}\to
 \nu\circ\iota^{-1}\circ\pi_k^{-1}$, weakly as $N\to\infty$
\cite[Exercises 30.3, 31.2]{Bauer}.

Let 
\[
{\cal T} = \bigcup_{k\in\nintegers} 
\{ f_{k,i}\circ\pi_k\circ\iota:\ \msp\times[0,1)\to\reals \;;\ 
f_{k,i}\in {\cal T}_k \},
\]
be the ${\cal T}$ in the assumption of \lemu{HDLproof0}. Since
$f_{k,i}$, $\pi$, $\iota$ are continuous, the functions in ${\cal T}$
are continuous. Note further that since $f_{k,i}$ is of bounded support,
the functions in ${\cal T}$ are uniformly continuous.
Since a countable union of countable sets is countable,
${\cal T}$ so defined is a countable set.
With this choice of ${\cal T}$, the assumption \eqnu{HDLproof0},
with a change in integration variable $z=\pi_k\circ\iota(\rho,y)$,
implies
\[ \arrb{l}\dsp
\limf{N} 
\int_{\nreals{k}\times[0,1)} f_{k,i}(z)\,
\nu_N\circ\iota^{-1}\circ\pi_k^{-1}(dz)=
\limf{N} 
\int_{\msp\times[0,1)}  f_{k,i}\circ\pi_k\circ\iota(\rho,y)\,
\nu_N(d\rho\times dy)
\\ \dsp {} =
\int_{\msp\times[0,1)}  f_{k,i}\circ\pi_k\circ\iota(\rho,y)\,
\nu(d\rho\times dy)
=
\int_{\nreals{k}\times[0,1)} f_{k,i}(z)\,
\nu\circ\iota^{-1}\circ\pi_k^{-1}(dz),
\arre \]
for all $k$ and $i$, which, as noted below \eqnu{rev0724-7},
implies $\nu_N\circ\iota^{-1}\circ\pi_k^{-1}\to
\nu\circ\iota^{-1}\circ\pi_k^{-1}$, weakly as $N\to\infty$,
for all $k$.
This implies that as measures on $\nreals{\infty}\times[0,1)$,
$\nu_N\circ\iota^{-1}\to\nu\circ\iota^{-1}$, weakly as $N\to\infty$
\cite[\S2 Example 2.4]{Billingsley}.
As noted in the paragraph between \eqnu{rev0724-4} and \eqnu{rev0724-6},
this further implies $\nu_N\to\nu$, weakly as $N\to\infty$.
\prfe
\remb
We could alternatively make use of separability of 
${\cal M}(\preals)$ directly to obtain a countable set ${\cal T}$,
following the discussion in 
\cite[\S 1, Remark 4.17, and remark after Corollary 9.3]{Kotani}.
\DDD\reme
\prfofb{\protect\lemu{Bauer}}
Let ${\cal T}$ be as in \lemu{HDLproof0}. If 
there exists, for each $n\in\nintegers$,
$\tilde{\Omega}_n\subset\Omega$ such that \eqnu{HDLproof0} holds
for $\omega\in\tilde{\Omega}_n$ and $\prb{\tilde{\Omega}_n}=1$ holds, then
$\Omega':=\bigcap_{n=1}^{\infty} \tilde{\Omega}_n$ satisfies
$\prb{\Omega'}=1$ and 
\eqnu{HDLproof0} holds for all $\omega\in\Omega'$ 
and $f_n\in {\cal T}$, which,
with \lemu{HDLproof0},  
implies \thmu{HDLinhomog}.

Let $d$ be the metric on $\msp\times[0,1)$ 
as in the proof of \lemu{HDLproof0}. 
Let $f_n\in {\cal T}$. Since $f_n$ is 
uniformly continuous, 
for any $\eps>0$ there exists $\delta>0$
such that for any $\rho_1,\rho_2\in\msp$ and $y_1,y_2\in [0,1)$, 
$d((\rho_1,y_1),(\rho_2,y_2))<\delta$ implies 
$|f_n(\rho_1,y_1)-f_n(\rho_2,y_2)|<\eps$.
Let $k$ be a positive integer greater than $1/{\delta}$ and put
\eqnb
\label{e:Bauerproof}
 f_{n,k}(\rho,y) =\sum_{l=0}^{k-1} f_{n}(\rho,{l}/{k})
 \chi_{[{l}/k,(l+1)/k)}(y),
\end{equation}
where $\chi_{[a,b)}(y)=1$ if $a\le y<b$ and $0$ otherwise.
Then for each $\rho\in \msp$ we have
\[ \sup_{y\in[0,1)} |f_n(\rho,y)-f_{n,k}(\rho,y)|<\eps. \]
Therefore,
$\lim_{k\to\infty} f_{n,k} =f_n$ uniformly on $\msp\times[0,1)$.
Noting that 
\[ \chi_{[l/k,(l+1)/k)}=
 \chi_{[{l}/k,1)}- \chi_{[(l+1)/k,1)}\,, \]
we see from \eqnu{Bauerproof} that
$f_{n,k}$ has an expression
\[ f_{n,k}(\rho,y)
=\sum_{l=0}^{k-1} g_{n,k,l}(\rho) \chi_{[{l}/k,1)}(y), \]
where $g_{n,k,l}:\ \msp\to\reals$ is bounded continuous.

Therefore, if \eqnu{HDLproof1} holds, then
using the definition \eqnu{muN} and the explicit formula \eqnu{HDLdistri} 
claimed in \thmu{HDLinhomog}, we see that there exists $\tilde{\Omega}_{n,k}$
satisfying $\prb{\tilde{\Omega}_{n,k}}=1$ and
\[
\limf{N} 
\int_{\msp\times[0,1)} f_{n,k}(\rho,y)\, \mu\ssN_t(d\rho\times dy)(\omega)
= \int_{\msp\times[0,1)} f_{n,k}(\rho,y)\, \mu_t(d\rho\times dy)
\]
ifl $\omega\in\tilde{\Omega}_{n,k}$\,.
Hence, $\tilde{\Omega}_n=\bigcap_{k=1}^{\infty}\tilde{\Omega}_{n,k}$
satisfies $\prb{\tilde{\Omega}_n}=1$ and \eqnu{HDLproof0} holds
for $\omega\in\tilde{\Omega}_n$\,.
\QED\prfofe

In view of \lemu{Bauer},
we fix $(y,t)$ and a bounded continuous function $g$,
in the remainder of this section.
Since $g$ is bounded, there exists a constant $M>0$ such that
\eqnb
\label{e:gM}
|g(\rho)|\le M, \ \rho\in\msp.
\end{equation}
Since the jump times $\{\tau\ssN_{i,1}\}$ are independent,
\prpu{yC} is proved in a straightforward way.
In contrast, $\{Y\ssN_i\}$ appearing in the left hand side of \eqnu{HDLproof1}
are dependent, and moreover, the 
non-linearity in \eqnu{Burgers} indicates that the dependence
cannot be neglected in the limit $N\to\infty$.
A strategy, inherited from the proof in \cite{HH071}, is to 
(i) choose a nice quantity defined as a sum of independent random variables
in such a way that the quantity converges to the right hand side of 
\eqnu{HDLproof1}, and (ii) show that the difference between the chosen 
quantity and the left hand side of \eqnu{HDLproof1}
can be shown to disappear in the limit, using the properties of the model.
We state these two steps explicitly in the following two Lemmas, 
respectively.
\lemb
\lema{HDLproof2} 
The following hold.
\itmb
\item
For $0\le y\le y_C(t)$,
\eqnb
\label{e:HDLproof21}
\frac1N\sum_{i=1}^Ng(\rho\ssN_i)\chrfcn{\nu\ssN_i((t-t_0(y,t), t])>0}
\ \to\ \int_{\msp} g(\rho)\,(1-e^{-\rho((t-t_0(y,t),t])})\, \Lambda(d\rho),
\end{equation}
almost surely as $N\to\infty$.

\item
For $y_C(t)\le y<1$,
\eqnb
\label{e:HDLproof22}
\frac1N\sum_{i=1}^Ng(\rho\ssN_i)
\chrfcn{(x\ssN_i-1)/N\ge \hat{y}(y,t),\ \tau\ssN_{i,1}>t}
\ \to\ 
\int_{\msp}g(\rho)\,e^{-\rho((0,t])}\,\mu_{0}(d\rho\times [\hat{y}(y,t),1)),
\end{equation}
almost surely as $N\to\infty$.
\DDD\itme
\leme

\lemb
\lema{HDLproof3} 
The following hold.
\itmb
\item
For $0\le y\le y_C(t)$,
\eqnb
\label{e:HDLproof31}
\frac1N\sum_{i=1}^N 
|\chrfcn{Y\ssN_i(t)< y}- \chrfcn{\nu\ssN_i((t-t_0(y,t), t])>0}|
\to 0,
\end{equation}
almost surely as $N\to\infty$.

\item
For $y_C(t)\le y<1$,
\eqnb
\label{e:HDLproof32}
\frac1N\sum_{i=1}^N 
|\chrfcn{Y\ssN_i(t)\ge y}
- \chrfcn{(x\ssN_i-1)/N \ge \hat{y}(y,t),\ \tau\ssN_{i,1}>t}|
\to 0,
\end{equation}
almost surely as $N\to\infty$.
\DDD\itme
\leme

\prfofb{\protect\eqnu{HDLproof1} assuming 
\protect\lemu{HDLproof2} and  \protect\lemu{HDLproof3}}
For the case $0\le y\le y_C(t)$, 
\eqnu{gM}, \eqnu{LambdaN}, \eqnu{HDLproof21}, and \eqnu{HDLproof31} imply
\[ \arrb{l}\dsp
\biggl|
\frac1N\sum_{i=1}^N g(\rho\ssN_i)\, \chrfcn{Y\ssN_i(t)\ge y} 
-\int_{\msp} g(\rho)\,e^{-\rho((t-t_0(y,t),t])}\, \Lambda(d\rho) 
\biggr|
\\ \dsp {}=
\biggl|
\frac1N\sum_{i=1}^N g(\rho\ssN_i)\, (1-\chrfcn{Y\ssN_i(t)< y})
-\int_{\msp} g(\rho)\,e^{-\rho((t-t_0(y,t),t])}\, \Lambda(d\rho) 
\biggr|
\\ \dsp {}=
\biggl|
\frac1N\sum_{i=1}^N g(\rho\ssN_i)\,
(-\chrfcn{Y\ssN_i(t)< y}+\chrfcn{\nu\ssN_i((t-t_0(y,t), t])>0})
\\ \dsp \phantom{{}=} +
\biggl(
\frac1N\sum_{i=1}^N g(\rho\ssN_i)-\int_{\msp} g(\rho)\, \Lambda(d\rho) 
\biggr)
\\ \dsp \phantom{{}=} +
\biggl(
-\frac1N\sum_{i=1}^N g(\rho\ssN_i)\,\chrfcn{\nu\ssN_i((t-t_0(y,t), t])>0}
+\int_{\msp} g(\rho)\,(1-e^{-\rho((t-t_0(y,t),t])})\, \Lambda(d\rho) 
\biggr)
\biggr|
\\ \dsp {} \le
M\, \frac1N\sum_{i=1}^N 
|\chrfcn{Y\ssN_i(t)< y}-\chrfcn{\nu\ssN_i((t-t_0(y,t), t])>0}|
\\ \dsp \phantom{{}=} +
\biggl|
\int_{\msp} g(\rho)\,\Lambda\ssN(d\rho)-\int_{\msp} g(\rho)\, \Lambda(d\rho)
\biggr|
\\ \dsp \phantom{{}=} +
\biggl|
\frac1N\sum_{i=1}^N g(\rho\ssN_i)\,\chrfcn{\nu\ssN_i((t-t_0(y,t), t])>0}
-\int_{\msp} g(\rho)\,(1-e^{-\rho((t-t_0(y,t),t])})\, \Lambda(d\rho) 
\biggr|
\\ \dsp
\ \to 0,\ \mbox{ a.s., as } N\to\infty,
\arre\]
which proves \eqnu{HDLproof1} for $0\le y\le y_C(t)$.

Similarly, for the case $y_C(t)\le y<1$,
\eqnu{gM}, \eqnu{HDLproof22}, and \eqnu{HDLproof32} imply
\[ \arrb{l}\dsp
\biggl|
\frac1N\sum_{i=1}^N g(\rho\ssN_i)\, \chrfcn{Y\ssN_i(t)\ge y} 
-\int_{\msp}g(\rho)\,e^{-\rho((0,t])}\,\mu_{0}(d\rho\times [\hat{y}(y,t),1))
\biggr|
\\ \dsp {}=
\biggl|
\frac1N\sum_{i=1}^N g(\rho\ssN_i)\, 
\biggl(
\chrfcn{Y\ssN_i(t)\ge y}
-\chrfcn{(x\ssN_i-1)/N \ge \hat{y}(y,t),\ \tau\ssN_{i,1}>t}
\biggr)
\\ \dsp \phantom{{}=} +
\biggl(
\frac1N\sum_{i=1}^N g(\rho\ssN_i)
 \chrfcn{(x\ssN_i-1)/N \ge \hat{y}(y,t),\, \tau\ssN_{i,1}>t}
\\ \dsp \phantom{{}=+(} 
-\int_{\msp}g(\rho) e^{-\rho((0,t])}\,\mu_{0}(d\rho\times [\hat{y}(y,t),1))
\biggr)
\biggr|
\\ \dsp {}\le
M\, \frac1N\sum_{i=1}^N 
|\chrfcn{Y\ssN_i(t)\ge y}
- \chrfcn{(x\ssN_i-1)/N \ge \hat{y}(y,t),\ \tau\ssN_{i,1}>t}|
\\ \dsp \phantom{{}=} +
\biggl|
\frac1N\sum_{i=1}^N g(\rho\ssN_i)
 \chrfcn{(x\ssN_i-1)/N \ge \hat{y}(y,t),\, \tau\ssN_{i,1}>t}
\\ \dsp \phantom{{}=+(} 
-\int_{\msp}g(\rho) e^{-\rho((0,t])}\,\mu_{0}(d\rho\times [\hat{y}(y,t),1))
\biggr|
\ \to 0,\ \mbox{ a.s., as } N\to\infty,
\arre\]
which proves \eqnu{HDLproof1} for $y_C(t)\le y<1$.
\QED
\prfofe

Before proving \lemu{HDLproof2} and  \lemu{HDLproof3},
we prepare a couple of random variables which converge as $N\to\infty$
to $y_A$ in \eqnu{yA} and $y_B$ in \eqnu{yB}.
The following \lemu{yANyBN} is used in the proof of \lemu{HDLproof3},
and the proof of \lemu{HDLproof2} is similar to that of \lemu{yANyBN}.
\lemb
\lema{yANyBN}
\itmb
\item
For $0\le t_0\le t$ define
\eqnb
\label{e:yAN}
Y\ssN_A(t_0,t)=\frac1N\sum_{i=1}^N \chrfcn{\nu\ssN_i((t-t_0,t])>0}.
\end{equation}
Then $Y\ssN_A(t_0,t)\to y_A(t_0,t)$, almost surely as $N\to\infty$.

\item
For $t\ge 0$ and $0\le y_0<1$ define
\eqnb
\label{e:yBN}
Y\ssN_B(y_0,t)=y_0+\frac1N\sum_{i;\ (x\ssN_i-1)/N\ge y_0}
\chrfcn{\tau\ssN_{i,1}\le t}.
\end{equation}
Then $Y\ssN_B(y_0,t)\to y_B(y_0,t)$, almost surely as $N\to\infty$.
\DDD

\itme
\leme
\prfb
As in the proof of \prpu{yC}, a strong law of large numbers implies,
almost surely as $N\to\infty$
\[ 
Y\ssN_A(t_0,t)-\EE{Y\ssN_A(t_0,t)}\to 0\ \mbox{ and }
Y\ssN_B(y_0,t)-\EE{Y\ssN_B(y_0,t)}\to 0.
\]
On the other hand, \eqnu{lambdaNstLambda} and \eqnu{lambdaN2lambda} imply
\[ \arrb{l}\dsp
\limf{N}\EE{Y\ssN_A(t_0,t)}=
\limf{N}\frac1N\sum_{i=1}^N (1-\prb{\nu\ssN_i((t-t_0,t])=0})
\\ \dsp {}
= \limf{N}\frac1N\sum_{i=1}^N (1-e^{-\rho\ssN_i((t-t_0,t])})
=1- \limf{N}\int_{\msp} e^{-\rho((t-t_0,t])} \Lambda\ssN(d\rho)
\\ \dsp {} 
=1- \limf{N}\int_0^{\infty} e^{-s} \lambda\ssN_{t-t_0,t}(ds)
=1- \int_{\msp} e^{-\rho((t-t_0,t])} \Lambda(d\rho)
= y_A(t_0,t).
\arre \]
Similarly,
\[ \arrb{l}\dsp
\limf{N}\EE{Y\ssN_B(y_0,t)}=
y_0+
\limf{N}\frac1N\sum_{i;\ (x\ssN_i-1)/N\ge y_0}\prb{\tau\ssN_{i,1}\le t}
\\ \dsp {}
=y_0+
\limf{N}\frac1N\sum_{i;\ (x\ssN_i-1)/N\ge y_0}
(1- e^{-\rho\ssN_i((0,t])})
\\ \dsp {}
=y_0+1- y_0
- \limf{N}\int_{\msp} e^{-\rho((0,t])} \mu\ssN_0(d\rho\times[y_0,1))
\\ \dsp {}
=1-\int_{\msp} e^{-\rho((0,t])} \mu_0(d\rho\times[y_0,1))
= y_B(y_0,t).
\arre \]
\prfe

\prfofb{\protect\lemu{HDLproof2}}
The proof is a repetition of the proof of \lemu{yANyBN},
by replacing $\chrfcn{\nu\ssN_i((t-t_0,t])>0}$
with $
g(\rho\ssN_i)\chrfcn{\nu\ssN_i((t-t_0(y,t), t])>0}$ for \eqnu{HDLproof21},
and $\chrfcn{\tau\ssN_{i,1}\le t}$
with $
g(\rho\ssN_i)\chrfcn{\tau\ssN_{i,1}>t}$ for \eqnu{HDLproof22}.
\QED
\prfofe

The proof of \eqnu{HDLproof1} now will be complete if we prove
\lemu{HDLproof3}, which is proved in a similar way as the corresponding
part in \cite{HH071}.
\prfofb{\protect\eqnu{HDLproof31} for $0\le y\le y_C(t)$}
Note that $y_A(t_0,t)$ of \eqnu{yA} is non-decreasing in $t_0$ and $t$,
with $y_A(0,t)=0$ and $y_A(t,t)=y_C(t)$,
and by assumption of the \thmu{HDLinhomog}, is continuous.
Hence 
\eqnb
\label{e:t0range}
y_A(t_0(y,t),t)=y,\ \ 0\le y\le y_C(t),\ t\ge0.
\end{equation}

\lemu{yANyBN} therefore implies that there exists $\Omega_A\subset\Omega$,
satisfying $\prb{\Omega_A}=1$, such that
\eqnb
\label{e:OmegaA}
\limf{N} Y\ssN_A(t_0(y,t),t)(\omega)=y,\ \ \omega\in\Omega_A\,.
\end{equation}
Fix $\omega\in \Omega_A$ arbitrarily.
The definition of the stochastic ranking process and \eqnu{yAN} imply
that $\nu\ssN_i((t-t_0(y,t), t])(\omega)>0$, if and only if
$Y\ssN_i(t)(\omega)$ is on the top side of $Y\ssN_A(t_0(y,t),t)(\omega)$;
$Y\ssN_i(t)(\omega)< Y\ssN_A(t_0(y,t),t)(\omega)$.
Therefore
\eqnb
\label{e:HDLproof311}\ \ \ \ 
\arrb{l}\dsp
\frac1N\sum_{i=1}^N 
|\chrfcn{Y\ssN_i(t)< y}(\omega)-\chrfcn{\nu\ssN_i((t-t_0(y,t), t])>0}(\omega)|
\\ \dsp {}
=
\frac1N\sum_{i=1}^N 
|\chrfcn{Y\ssN_i(t)< y}(\omega)
- \chrfcn{Y\ssN_i(t)< Y\ssN_A(t_0(y,t),t)}(\omega)|.
\arre
\end{equation}
Note that the definition of $Y\ssN_i(t)$ in \eqnu{YN} implies
that it takes values in $\{{k}/N \;;\ k=0,1,\ldots,N-1\}$.
Hence \eqnu{OmegaA} implies
\eqnb
\label{e:HDLproof312}
\arrb{l}\dsp
\limsup_{N\to\infty} \frac1N\sum_{i=1}^N 
|\chrfcn{Y\ssN_i(t)< y}(\omega)
- \chrfcn{Y\ssN_i(t)< Y\ssN_A(t_0(y,t),t)}(\omega)|
\\ \dsp {}
\le
\limf{N} \frac1N\times (N |Y\ssN_A(t_0(y,t),t)(\omega)-y|+1)
= 0.
\arre
\end{equation}
The relations \eqnu{HDLproof311} and \eqnu{HDLproof312}
imply \eqnu{HDLproof31}.
\QED\prfofe

\prfofb{\protect\eqnu{HDLproof32} for $y_C(t)\le y<1$}
 $y_B(y,t)$ of \eqnu{yB} is non-decreasing in $y$ and $t$,
with $y_B(0,t)=y_C(t)$ and $y_B(1-,t)=1-0$,
and by assumption of the \thmu{HDLinhomog}, is continuous.
Hence 
\eqnb
\label{e:hatyrange}
y_B(\hat{y}(y,t),t)=y,\ \ y_C(t)\le y<1,\ t\ge0.
\end{equation}
\lemu{yANyBN} therefore implies that there exists $\Omega_B\subset\Omega$,
satisfying $\prb{\Omega_B}=1$, such that
\eqnb
\label{e:OmegaB}
\limf{N} Y\ssN_B(\hat{y}(y,t),t)(\omega)=y,\ \ \omega\in\Omega_B\,.
\end{equation}
Fix $\omega\in \Omega_B$ arbitrarily.
The definition of the stochastic ranking process and \eqnu{yBN} imply
that 
$(x\ssN_i-1)/N \ge \hat{y}(y,t)$ and $\tau\ssN_{i,1}(\omega)>t$
hold together, if and only if
$Y\ssN_i(t)(\omega)$ is on the tail side of $Y\ssN_B(\hat{y}(y,t),t)(\omega)$;
$Y\ssN_i(t)(\omega)\ge Y\ssN_B(t_0(y,t),t)(\omega)$.
Therefore
\eqnb
\label{e:HDLproof321}
\arrb{l}\dsp
\frac1N\sum_{i=1}^N 
|\chrfcn{Y\ssN_i(t)\ge y}(\omega)
- \chrfcn{(x\ssN_i-1)/N \ge \hat{y}(y,t),\ \tau\ssN_{i,1}>t}(\omega)|
\\ \dsp
=
\frac1N\sum_{i=1}^N 
|\chrfcn{Y\ssN_i(t)\ge y}(\omega)
- \chrfcn{Y\ssN_i(t)\ge Y\ssN_B(\hat{y}(y,t),t)}(\omega)|.
\arre
\end{equation}
As in the proof of \eqnu{HDLproof31}, $Y\ssN_i(t)$ 
takes values in $\{{k}/N \;;\ k=0,1,\ldots,N-1\}$,
which implies, with \eqnu{OmegaB},
\eqnb
\label{e:HDLproof322}
\arrb{l}\dsp
\limsup_{N\to\infty} \frac1N\sum_{i=1}^N 
|\chrfcn{Y\ssN_i(t)\ge y}(\omega)
- \chrfcn{(x\ssN_i-1)/N \ge \hat{y}(y,t),\ \tau\ssN_{i,1}>t}(\omega)|
\\ \dsp
\le
\limf{N} \frac1N\times (N |Y\ssN_B(\hat{y}(y,t),t)(\omega)-y|+1)
= 0.
\arre\end{equation}
The relations \eqnu{HDLproof321} and \eqnu{HDLproof322}
imply \eqnu{HDLproof32}.
\QED\prfofe
This completes the proof of \lemu{HDLproof3}, hence of \thmu{HDLinhomog}.

\section{Proof of \protect\thmu{Burgers}.}
\seca{Burgers}
To prove \thmu{Burgers}, we apply a standard method of characteristic curves.

First, assume $0\le y\le y_C(t)=y_A(t,t)$.
Let $t_1\ge 0$, and consider an ordinary differential equation 
for a characteristic curve intersecting $(0,t_1)$, defined by
\eqnb
\label{e:BurgersproofA1}
\arrb{l}\dsp
\diff{y}{t}(t)=\sum_{\beta=1}^k w_{\beta}(t)\, U_{\beta}(y(t),t),
\ \ \alpha=1,2,\ldots,k,\ t\ge t_1\,,
\\ \dsp
y(t_1)=0.
\arre
\end{equation}
Put
\eqnb
\label{e:BurgersproofA2}
\varphi_{\alpha}(t)=U_{\alpha}(y(t),t),
\ \ \alpha=1,2,\ldots,k,\ t\ge t_1\,.
\end{equation}
Then \eqnu{BurgersproofA2}, \eqnu{Burgers}, and \eqnu{BurgersproofA1} imply
\eqnb
\label{e:BurgersproofA3}
\diff{\varphi_{\alpha}}{t}(t)
=-w_{\alpha}(t) U_{\alpha}(y(t),t)
=-w_{\alpha}(t) \varphi_{\alpha}(t),
\end{equation}
which, with $y(t_1)=0$ in \eqnu{BurgersproofA1}, has a unique solution 
\eqnb
\label{e:BurgersproofA4}
\varphi_{\alpha}(t)
=U_{\alpha}(0,t_1)\,\exp( -\int_{t_1}^t w_{\alpha}(u)du )
=r_{\alpha}\,\exp( -\int_{t_1}^t w_{\alpha}(u)du ),
\end{equation}
where we also used \eqnu{Burgersboundary}.
Substituting \eqnu{BurgersproofA2} and \eqnu{BurgersproofA4} in
\eqnu{BurgersproofA1}, we have
\[
\diff{y}{t}(t)=\sum_{\beta=1}^k 
w_{\beta}(t)\, r_{\beta}\,\exp( -\int_{t_1}^t w_{\beta}(u)du ),
\]
which, with $y(t_1)=0$, has a unique solution
\eqnb
\label{e:BurgersproofA5}
y(t)=\sum_{\beta=1}^k 
r_{\beta}\,(1-\exp( -\int_{t_1}^t w_{\beta}(u)du ))
= y_A(t-t_1,t).
\end{equation}
where we also used $\sum_{\beta=1}^k r_{\beta}=1$ in \eqnu{Burgersassump}
and \eqnu{BurgersyA} with \eqnu{Burgersrho}, in the last equality.
The assumptions for $w_{\alpha}$ in \thmu{Burgers} imply
that $y_A(t_0,t)$ is strictly increasing and differentiable in $t_0$\,,
satisfying $y_A(0,t)=0$ and $y_A(t,t)=y_C(t)$.
Hence there exists a unique, strictly increasing, 
differentiable inverse function $t_0=t_0(y,t)$, taking values in $[0,t]$,
satisfying 
\[ y_A(t_0(y,t),t)=y,\ 0\le y\le y_C(t),\ t\ge0.  \]
This, with \eqnu{BurgersproofA2}, \eqnu{BurgersproofA4}, and
\eqnu{BurgersproofA5}, implies
\[
U_{\alpha}(y,t)
=r_{\alpha}\,\exp( -\int_{t-t_0(y,t)}^t w_{\alpha}(u)du ),
\]
which proves \eqnu{HDLdistridiscrete} for $0\le y\le y_C(t)$.

Next, assume $y_C(t)=y_B(0,t)\le y<1$.
Let $0\le y_0<1$, and consider an ordinary differential equation 
for a characteristic curve intersecting $(y_0,0)$, defined by
\eqnb
\label{e:BurgersproofB1}
\arrb{l}\dsp
\diff{y}{t}(t)=\sum_{\beta=1}^k w_{\beta}(t)\, U_{\beta}(y(t),t),
\ \ \alpha=1,2,\ldots,k,\ t\ge 0\,,
\\ \dsp
y(0)=y_0\,.
\arre
\end{equation}
Put
\eqnb
\label{e:BurgersproofB2}
\varphi_{\alpha}(t)=U_{\alpha}(y(t),t),
\ \ \alpha=1,2,\ldots,k,\ t\ge t_1\,.
\end{equation}
Then \eqnu{BurgersproofB2}, \eqnu{Burgers}, and \eqnu{BurgersproofB1} imply,
exactly as for the case $y\le y_C(t)$,
\eqnb
\label{e:BurgersproofB3}
\diff{\varphi_{\alpha}}{t}(t)=-w_{\alpha}(t) \varphi_{\alpha}(t),
\end{equation}
which, with $y(0)=y_0$, has a unique solution 
\eqnb
\label{e:BurgersproofB4}
\varphi_{\alpha}(t)=u_{\alpha}(y_0)\,\exp( -\int_{0}^t w_{\alpha}(u)du ),
\end{equation}
where we also used \eqnu{Burgersinitial}.
Substituting \eqnu{BurgersproofB2} and \eqnu{BurgersproofB4} in
\eqnu{BurgersproofB1}, we have another differential equation for $y(t)$,
which, with $y(0)=y_0$, has a unique solution
\eqnb
\label{e:BurgersproofB5}
y(t)
= y_B(y_0,t),
\end{equation}
where we used
$\sum_{\beta=1}^k u_{\beta}(y) = 1-y$ in \eqnu{Burgersassump}
and \eqnu{BurgersyB} with \eqnu{Burgersrho}.
The assumptions for $u_{\alpha}$ in \thmu{Burgers} imply
that $y_B(y,t)$ is strictly increasing and differentiable in $y$\,,
satisfying $y_B(0,t)=y_C(t)$ and $y_B(1-,t)=1-$.
Hence there exists a unique, strictly increasing, 
differentiable inverse function $\hat{y}(y,t)$, taking values in $[0,1)$,
satisfying 
\[ y_B(\hat{y}(y,t),t)=y,\ y_C(t)\le y< 1,\ t\ge0.  \]
As in the proof for $y\le y_C(t)$,
this, with \eqnu{BurgersproofB2}, \eqnu{BurgersproofB4}, and
\eqnu{BurgersproofB5}, implies \eqnu{HDLdistridiscrete} for $ y_C(t)\le y<1$.

This completes a proof of \thmu{Burgers}.\QED

\section{Scaling limit results uniform in time.}
\seca{fnctHDL}

Let $T>0$ and
\eqnb
\label{e:equicontiintensity}
{\cal I}
=\{r\ssN_i:\ [0,T]\to \preals \;;\ i=1,2,\ldots,N,\ N\in\nintegers \}
\end{equation}
be a set of continuous functions on $[0,T]$ defined by
$r\ssN_i(t)=\rho\ssN_i((0,t])$, $t\ge0$.
Note that since we assumed in the beginning that $\rho\ssN_i$ is continuous,
$r\ssN_i$ is continuous. In this section, we prove the following.
\thmb
\thma{fnctHDLinhomog}
Let $T>0$.
In addition to the assumptions in  \prpu{yC}, assume that 
a set of continuous functions ${\cal I}$ defined by \eqnu{equicontiintensity}
is uniformly equicontinuous; namely, 
\eqnb
\label{e:fnctHDLinhomog1}
\lim_{\delta\downarrow0} \sup_{r\in {\cal I}} 
 \sup_{s,t\in[0,T];\ |s-t|\le\delta} |r(s)-r(t)| = 0. 
\end{equation}
Then, $Y\ssN_C$ of \eqnu{YCN} converges almost surely to $y_C$ of
\eqnu{yC} as $N\to\infty$, as a sequence in the space of continuous
functions on $[0,T]$ with supremum norm\/$:$
\eqnb
\label{e:fnctyC}
\mathrm{P}\left[\;
\limf{N} \sup_{t\in[0,T]} |Y\ssN_C(t)-y_C(t)| =0
\;\right]=1.
\end{equation}

Assume next that all
the assumptions of \thmu{HDLinhomog} and \eqnu{fnctHDLinhomog1}
hold. Assume also that a set of functions
\[
{\cal J}=\{r:\ [0,T]\to\preals\;;\ r(t)=\rho((0,t]),\ t\in[0,T],
\ \rho\in \suppt \Lambda\}
\]
is uniformly equicontinuous, and
that for $y_A$ of \eqnu{yA} and 
$y_B$ of \eqnu{yB}, $y_A(t-t_1,t)$ and $y_B(y,t)$ are equicontinuous in 
$(t_1,t)$ and $(y,t)$, respectively. Then,
$\mu\ssN_{\cdot}$ of \eqnu{muN} converges almost surely to $\mu_{\cdot}$
of \eqnu{HDLdistri} as $N\to\infty$, as a sequence in the space of probability
measure valued functions $\mu_{\cdot}:\ t\mapsto\mu_t$ with supremum norm.
\DDD\thme
\prfb
First we assume that the assumptions of \prpu{yC} and 
\eqnu{fnctHDLinhomog1} hold.
Note that \eqnu{PRM} implies that, for $i=1,2,\ldots,N$, $N=1,2,\ldots$,
\[
 \nu\ssN_i((0,t])-r\ssN_i(t)
= \nu\ssN_i((0,t])-\rho\ssN_i((0,t]),\ \ t\ge0,
\]
is a martingale up to fixed time $T$.
Note also that \eqnu{PoissonUnitMeasRep} implies
\eqnb
\label{e:YCintensitystoppingrep}
 \nu\ssN_i((0,t\wedge \tau\ssN_{i,1}])= \chrfcn{\tau\ssN_{i,1}\le t}\,.
\end{equation}
Hence 
\eqnb
\label{e:fnctHDLinhomogprf4}
W\ssN_i(t):=
 \chrfcn{\tau\ssN_{i,1}\le t} -r\ssN_i(t\wedge\tau\ssN_{i,1}),
\ t\in [0,T],
\end{equation}
is a bounded martingale. This with \eqnu{YCN} further implies that
\[ Y_C^N(t) - \frac{1}{N} \sum_{i=1}^N r\ssN_i(t\wedge\tau\ssN_{i,1}) 
= \frac{1}{N} \sum_{i=1}^N 
(\chrfcn{\tau\ssN_{i,1}\le t} -r\ssN_i(t\wedge\tau\ssN_{i,1}))
=\frac1N\sum_{i=1}^N W\ssN_i(t)
 \]
is also a bounded martingale.
Using Doob's inequality, 
independence of $\{\tau\ssN_{i,1}\;;\ i=1,2,\ldots,N\}$,
and $|W\ssN_i(T)|\le1$, we have
\[ 
\mathrm{E}\left[\;
\sup_{0 \le t \le T} \biggl(\frac{1}{N} \sum_{i=1}^N W\ssN_i(t)\biggr)^4
\;\right]
\le \frac{4^4}{3^4}\, 
\mathrm{E}\left[\;
\biggl(\frac{1}{N} \sum_{i=1}^N W\ssN_i(T)\biggl)^4
\;\right]
\le \frac{4^4}{3^3N^2}\,. \]
With an argument similar to that in the proof of \prpu{yC},
\eqnb
\label{e:fnctHDLinhomogprf1}
\sup_{0 \le t \le T} \biggl|\frac{1}{N} \sum_{i=1}^N W\ssN_i(t)\biggr|
\to 0,\ \mbox{ a.e., as }\ N\to\infty.
\end{equation}

On the other hand, for each $0\le t\le T$, as in the proof of \prpu{yC},
independence and boundedness of
$r\ssN_i(t\wedge\tau\ssN_{i,1})$, $i=1,2,\ldots,N$, imply
\eqnb
\label{e:fnctHDLinhomogprf6}
\frac1N\sum_{i=1}^Nr\ssN_i(t\wedge\tau\ssN_{i,1})
-\frac1N\sum_{i=1}^N\EE{r\ssN_i(t\wedge\tau\ssN_{i,1})}
\to 0,\ \mbox{ a.e., as }\ N\to\infty,
\end{equation}
and
\[
\EE{r\ssN_i(t\wedge\tau\ssN_{i,1})}=\EE{\chrfcn{\tau\ssN_{i,1}\le t}}
-\EE{W\ssN_i(t)}=\EE{\chrfcn{\tau\ssN_{i,1}\le t}}
\]
implies
\eqnb
\label{e:fnctHDLinhomogprf5}
\limf{N} \frac1N\sum_{i=1}^N\EE{r\ssN_i(t\wedge\tau\ssN_{i,1})}
=\limf{N} \EE{Y\ssN_C(t)} =y_C(t).
\end{equation}
Since $Y\ssN_C(t)$ is non-decreasing in $t$, and
$y_C(t)$ is its pointwise limit, it is also non-decreasing.
As in the case of \coru{yC}, \eqnu{fnctHDLinhomogprf6} and
\eqnu{fnctHDLinhomogprf5} imply that, with probability one,
\eqnb
\label{e:fnctHDLinhomogprf2}
\frac1N\sum_{i=1}^N r\ssN_i(t\wedge\tau\ssN_{i,1})\to y_C(t),
\ t\in\rationals\cap[0,T],\ \mbox{ as }\ N\to\infty.
\end{equation}
Since 
$\{r\ssN_i\}$ is equicontinuous, 
\eqnu{fnctHDLinhomogprf2} implies that $y_C$ is continuous on rationals,
and the monotonicity of $y_C$ proves that it is continuous on $[0,T]$.

By assumption of equicontinuity and
the convergence \eqnu{fnctHDLinhomogprf2} on a dense subset of $[0,T]$,
it follows that the convergence is uniform:
\eqnb
\label{e:fnctHDLinhomogprf3}
\sup_{t\in [0,T]}
\biggl|\frac{1}{N} \sum_{i=1}^N r\ssN_i(t\wedge \tau\ssN_{i,1})- y_C(t)
\biggr| \to 0,\ \mbox{ a.e., as }\ N\to\infty.
\end{equation}
The equations \eqnu{YCN}, \eqnu{fnctHDLinhomogprf4}, 
\eqnu{fnctHDLinhomogprf1}, \eqnu{fnctHDLinhomogprf3} prove \eqnu{fnctyC}.

In the remainder of this section, we assume that the assumptions of
\thmu{fnctHDLinhomog} hold.
To prove uniform convergence of $\mu\ssN_t$, we first prepare $t$-uniform 
version of \lemu{Bauer}.
\lemb
\lema{fnctBauer}
If, for each $y\in[0,1)$ and for each
bounded continuous function $g:\ \msp\to \reals$,
there exists $\tilde{\Omega}$ with $\prb{\tilde{\Omega}}=1$ such that,
for each $\omega\in\tilde{\Omega}$,
\eqnb
\label{e:fnctHDLproof0}
\limf{N} \sup_{t\in[0,T]}
\biggl|\int_{\msp} g(\rho)\,\mu\ssN_t(d\rho\times [y,1))(\omega)
-\int_{\msp} g(\rho)\,\mu_t(d\rho\times [y,1))\biggr|=0,
\end{equation}
then $\mu\ssN_t$ converges to $\mu_t$ uniformly in $t\in[0,T]$
as $N\to\infty$, almost surely.
\DDD\leme
\prfb
Let ${\cal T}=\{f_n\;;\ n\in\nintegers\}$ be as in the proof of \lemu{Bauer},
and for probability measures $\mu$ and $\nu$ on $\msp\times[0,1)$, put
\[ \arrb{l} \dsp
\pi(\mu,\nu)
\\ \dsp {}
:=\sum_{n=1}^{\infty} 2^{-n} \biggl(\biggl|
\int_{\msp\times[0,1)} f_n(\rho,y)\,\mu(d\rho\times dy) -
\int_{\msp\times[0,1)} f_n(\rho,y)\,\nu(d\rho\times dy) \biggr|\wedge 1
\biggr).
\arre \]
Then $\pi$ is a metric on the space of
probability measures on $\msp\times[0,1)$,
and the convergence with respect to $\pi$ is equivalent to
convergence \eqnu{HDLproof0} for each $f_n\in{\cal T}$.
Hence, as noted just below \eqnu{HDLproof0}, it is equivalent
to weak convergence of the probability measures on $\msp\times[0,1)$.

Now assume that \eqnu{fnctHDLproof0} holds.
Then following the arguments of the proof of \lemu{Bauer},
replacing \eqnu{HDLproof1} by \eqnu{fnctHDLproof0},
we see that there exists $\Omega'\subset\Omega$ such that $\prb{\Omega'}=1$
and 
\[
\limf{N} \sup_{t\in[0,T]} \biggl|
\int_{\msp\times[0,1)} f_n(\rho,y)\mu\ssN_t(d\rho\times dy)(\omega) -
\int_{\msp\times[0,1)} f_n(\rho,y)\mu_t(d\rho\times dy) \biggr|=0,
\]
for all $n\in\nintegers$ and $\omega\in\Omega'$.
Therefore,
\[ \limf{N} \sup_{t\in[0,T]}\pi(\mu\ssN_t(\omega),\mu_t)=0,
\ \ \omega\in\Omega', \]
which, by the equivalence of convergence in $\pi$ and the convergence
in the weak topology of the space of probability measures on
$\msp\times[0,1)$, implies the almost sure uniform convergence in $t\in[0,T]$,
of $\mu\ssN_t$ to $\mu_t$.
\prfe

In view of \lemu{fnctBauer},
we fix $y\in [0,1)$ and a bounded continuous function $g$,
in the remainder of this section. Note that \eqnu{gM} holds.
The assumption $\Lambda\ssN\to \Lambda$ in \eqnu{LambdaN} further implies
that for any $K>0$ there exists a positive integer $N_0$ such that, for
$N>N_0$\,,
\eqnb
\label{e:fnctHDLproof20}
\biggl|
\frac1N\sum_{i=1}^N g(\rho\ssN_i)-\int_{\msp} g(\rho)\, \Lambda(d\rho)
\biggr|
=
\biggl|
\int_{\msp} g(\rho)\,\Lambda\ssN(d\rho)-\int_{\msp} g(\rho)\, \Lambda(d\rho)
\biggr|
< \frac{M}K\,.
\end{equation}

The following Lemma corresponds to \lemu{yANyBN}.
\lemb
\lema{fnctyANyBN}
For each $t_1\in [0,T]$,
$Y\ssN_A$ of \eqnu{yAN} and $y_A$ of \eqnu{yA} satisfy
\eqnb
\label{e:yAno1}
 \sup_{t\in[t_1,T]} |Y\ssN_A(t-t_1,t)-y_A(t-t_1,t)| 
\to 0,\ \mbox{ a.e., as }\ N\to\infty, 
\end{equation}
and for each $y_0\in[0,1)$,
$Y\ssN_B$ of \eqnu{yBN} and $y_B$ of \eqnu{yB} satisfy
\eqnb
\label{e:yAno2}
 \sup_{t\in[0,T]} |Y\ssN_B(y_0,t)-y_B(y_0,t)| 
\to 0,\ \mbox{ a.e., as }\ N\to\infty. 
\end{equation}
\DDD\leme
\prfb
Define, for  $i=1,2,\ldots,N$, $N=1,2,\ldots$, 
\[\tilde{\tau} \ssN_{i}=  \tau\ssN_{i,k_i} ,\]
where $k_i := \inf \{j \;;\ \tau\ssN_{i,j} >t_1\}$.
Then just as in the proof of \eqnu{fnctyC}, we see  that 
\[V\ssN_i(t):=
 \chrfcn{\tilde{\tau}\ssN_{i}\le t}
 -\rho \ssN_i((t_1, t\wedge \tilde{\tau}\ssN_{i}]),
\ t\in [t_1,T]
\]
and, with \eqnu{yAN}, accordingly, 
\eqnb  \label{e:yAno3} Y_A^N(t-t_1,t ) - \frac{1}{N} \sum_{i=1}^N \rho 
\ssN_i((t_1, t\wedge \tilde{\tau}\ssN_{i}])
=\frac1N\sum_{i=1}^N V\ssN_i(t),
\ t\in[t_1,T] \end{equation}
are  bounded martingales, and we have
\eqnb\label{e:yAno4}
\sup_{t_1 \le t \le T} \biggl|  \frac1N\sum_{i=1}^N V\ssN_i(t) 
\biggr|
\to 0,\ \ t\in\rationals\cap[0,T],\ \ \mbox{ a.e., as }\ N\to\infty.
\end{equation}
On the other hand, we have with probability one, 
\eqnb\label{e:yAno5}
 \frac1N\sum_{i=1}^N \rho \ssN_i((t_1 , t\wedge\tilde{\tau}\ssN_i ])
=y_A(t-t_1,t),
\ t\in\rationals\cap[0,T],\ \mbox{ as }\ N\to\infty.
\end{equation}

By assumptions of  equicontinuity and
the convergence \eqnu{yAno5} on a dense subset of $[0,T]$,
it follows that the convergence is uniform:
\eqnb
\label{e:yAno6}
\sup_{t\in [0,T]}
\biggl|\frac{1}{N} \sum_{i=1}^N \rho \ssN_i((t_1, t\wedge \tau\ssN_i])
- y_A(t-t_1,t)
\biggr| \to 0,\ \mbox{ a.e., as }\ N\to\infty.
\end{equation}
The equations \eqnu{yAno3}, \eqnu{yAno4} and \eqnu{yAno6} 
lead to \eqnu{yAno1}. 

A proof of \eqnu{yAno2} goes in exact
correspondence with that of \eqnu{fnctyC}, if we directly use
the assumption of continuity of $y_B$ in place of monotonicity of $y_C$\,.
\prfe
\corb
\cora{fnctHDLproof2}
For each $t_1\in [0,T]$,
\eqnb
\label{e:fnctHDLproof21}
\sup_{t\in[t_1,T]} \biggl|
\frac1N\sum_{i=1}^Ng(\rho\ssN_i)\chrfcn{\nu\ssN_i((t_1, t])>0}
- \int_{\msp} g(\rho)\,(1-e^{-\rho((t_1,t])})\, \Lambda(d\rho)
\biggr|\to 0,
\end{equation}
almost surely as $N\to\infty$, 
and for each $y_0\in[0,1)$,
\eqnb
\label{e:fnctHDLproof22}
\sup_{t\in[0,T]} \biggl|
\frac1N\sum_{i=1}^Ng(\rho\ssN_i)
\chrfcn{(x\ssN_i-1)/N\ge y_0\,,\ \tau\ssN_{i,1}>t}-
\int_{\msp}g(\rho)\,e^{-\rho((0,t])}\,\mu_{0}(d\rho\times [y_0,1))
\biggr|\to 0,
\end{equation}
almost surely as $N\to\infty$.
\DDD\core
\prfb
This is proved as in the proof of \lemu{fnctyANyBN},
if one notes \eqnu{gM}.
\prfe

Fix a positive integer $K$ arbitrarily.
By the assumptions of \thmu{fnctHDLinhomog} of uniform equicontinuity of
${\cal J}$,  $y_A$ and $y_B$,
and noting that $\mu_0(\msp\times[y,1))=1-y$,
there exist a positive integer $L$ and sequences
$0=t_{1,0}<t_{1,1}<\cdots<t_{1,L}=T$ and 
$0=y_{0,0}<y_{0,1}<\cdots<y_{0,L}=1$ such that
\itmb
\item
for $j=0,1,2,\ldots,M$ and $s\in [t_{1,j-1}\,,\; t_{1,j+1}]$,
\eqnb
\label{e:fnctHDLproof29000}
\int_{\msp} | e^{-\rho((t_{1,j},t])}- e^{-\rho(s,t])} | 
\, \Lambda(d\rho)<\frac1K\,,
\end{equation}
where, for convenience we put $t_{1,j}=0$ if $j\le 0$, and $t_{1,j}=T$ if 
$j\ge L$, and 
also for $j=0,1,2,\ldots,M$ and $z\in [y_{0,j-1}\,,\; y_{0,j+1}]$,
\eqnb
\label{e:fnctHDLproof29001}
\arrb{l}\dsp
\left| \int_{\msp} g(\rho) e^{-\rho((0,t])} \mu_0(d\rho\times[y_{0,j},1))
- \int_{\msp} g(\rho) e^{-\rho((0,t])} \mu_0(d\rho\times[z,1)) \right|
\\ \dsp {} = 
\left| \int_{\msp} g(\rho) e^{-\rho((0,t])}
 \mu_0(d\rho\times[y_{0,j}\wedge z,y_{0,j}\vee z))\right|
<\frac{M}K\,,
\arre
\end{equation}
where, we put $y_{0,j}=0$ if $j\le 0$, and $y_{0,j}=1$ if $j\ge L$,

\item
the sequences of functions
$y_{A,j}(t)=y_A((t-t_{1,j})\vee0,t)$, $j=0,1,2,\ldots,L$,
which is decreasing in $j$,
and
$y_{B,j}(t)=y_B(y_{0,j},t)$, $j=0,1,2,\ldots,L-1$,
which is increasing in $j$,
satisfy
\eqnb
\label{e:fnctHDLproof23}
0\le y_{A,j}(t)-y_{A,j+1}(t)< \frac1K\,,\ j=0,1,2,\ldots,L-1,\ t\in[0,T],
\end{equation}
and
\eqnb
\label{e:fnctHDLproof24}
0\le y_{B,j+1}(t)-y_{B,j}(t)< \frac1K\,,\ j=0,1,2,\ldots,L-2,\ t\in[0,T].
\end{equation}
\itme

\lemu{fnctyANyBN} and \coru{fnctHDLproof2} imply that there exists 
$\tilde{\Omega}_K\subset\Omega$, satisfying $\prb{\tilde{\Omega}_K}=1$,
such that for all $\omega\in\tilde{\Omega}_K$ there exists an integer
$N_0=N_0(\omega)$ such that if $N>N_0$ then
\eqnb
\label{e:fnctHDLproof27}
| Y\ssN_A(t-t_{1,j},t)(\omega)-y_{A,j}(t) |<\frac1K\,,
\ t\in [t_{1,j},T],\ j=0,1,\ldots,L,
\end{equation}

\eqnb
\label{e:fnctHDLproof28}
| Y\ssN_B(y_{0,j},t)(\omega)-y_{B,j}(t) |<\frac1K\,,
\ j=0,1,\ldots,L,\ t\in [0,T],
\end{equation}

\eqnb
\label{e:fnctHDLproof25}
\arrb{l} \dsp
\biggl|
\frac1N\sum_{i=1}^Ng(\rho\ssN_i)\chrfcn{\nu\ssN_i((t_{1,j}, t])>0}(\omega)
- \int_{\msp} g(\rho)\,(1-e^{-\rho((t_{1,j},t])})\, \Lambda(d\rho)
\biggr|<\frac{M}K\,,\ \ 
\\ \dsp
t\in [t_{1,j},T],\ j=0,1,\ldots,L,
\arre
\end{equation}
and
\eqnb
\label{e:fnctHDLproof26}
\arrb{l} \dsp
\biggl|
\frac1N\sum_{i=1}^Ng(\rho\ssN_i)
\chrfcn{ (x\ssN_i-1)/N\ge y_{0,j}\,,\ \tau\ssN_{i,1}>t}(\omega)-
\int_{\msp}g(\rho)\,e^{-\rho((0,t])}\,\mu_{0}(d\rho\times [y_{0,j},1))
\biggr|
\\ \dsp {} 
<\frac{M}K\,,\ \ 
j=0,1,\ldots,L-1,\ t\in [0,T].
\arre
\end{equation}

Now, we shall consider the case $y_C(t)\ge y$ and the case $y_C(t)\le y$
separately.
First, let $y_C(t)\ge y$, and let $j=j(t)$ be the integer such that
\eqnb
\label{e:fnctHDLproof290}
y_{A,j}(t)\le y<y_{A,j-1}(t).
\end{equation}
Note that $y_C(t)\ge y$ implies $y=y_A(t_0(y,t),t)$ (see \eqnu{t0range}),
with which
$y_{A,0}(t)=y_A(t,t)=y_C(t)$, $y_{A,L}(t)=y_A(0,t)=0$,
and monotonicity of $y_A(t_0,t)$ with respect to $t_0$ imply that such
an integer $j=j(t)$ exists if $y_C(t)\ge y$.
Since $y_A(t_0,t)$ is increasing in $t_0$\,, \eqnu{fnctHDLproof290}
also implies 
\eqnb
\label{e:fnctHDLproof2900}
t_{1,j-1}<t-t_0(y,t)\le t_{1,j}.
\end{equation}
Since \eqnu{fnctHDLproof23} implies
\[ 0\le y-y_{A,j}(t) \le y_{A,j-1}(t)-y_{A,j}(t)<\frac1K\,, \]
with \eqnu{fnctHDLproof27} and a similar argument as for \eqnu{HDLproof312},
we have
\eqnb
\label{e:fnctHDLproof291}
\arrb{l}\dsp
 \biggl| \frac1N\sum_{i=1}^N g(\rho\ssN_i)\, 
(\chrfcn{Y\ssN_i(t)< y} - \chrfcn{Y\ssN_i(t)<Y\ssN_A(t-t_{1,j},t)})
(\omega) \biggr|
\\ \dsp {}
\le \frac1N\sum_{i=1}^N |g(\rho\ssN_i)|\,  \biggl| 
\chrfcn{Y\ssN_i(t)< y} - \chrfcn{Y\ssN_i(t)<y_{A,j}(t)} \biggr|(\omega)
\\ \dsp \phantom{\le}
+  \frac1N\sum_{i=1}^N |g(\rho\ssN_i)|\, \biggl|
\chrfcn{Y\ssN_i(t)<y_{A,j}(t)} - \chrfcn{Y\ssN_i(t)<Y\ssN_A(t-t_{1,j},t)}
 \biggr|(\omega)
\\ \dsp {}
 \le M (y-y_{A,j}(t))+ M|Y\ssN_A(t-t_{1,j},t)(\omega)-y_{A,j}(t) |
<\frac{2M}K\,.
\arre
\end{equation}
Note also that, as in the argument for \eqnu{HDLproof311},
\eqnb
\label{e:fnctHDLproof292}
\chrfcn{\nu\ssN_i((t_{1,j}, t])>0}=\chrfcn{Y\ssN_i(t)<Y\ssN_A(t-t_{1,j},t)}.
\end{equation}

Adding up \eqnu{fnctHDLproof20}, \eqnu{fnctHDLproof291},
\eqnu{fnctHDLproof25} and \eqnu{fnctHDLproof29000},
and using  \eqnu{fnctHDLproof292} and triangular inequality, we arrive at
\eqnb
\label{e:fnctHDLproof1}
\arrb{l} \dsp
\sup_{t\in[0,T];\ y_C(t)\ge y}\biggl|
\frac1N\sum_{i=1}^N g(\rho\ssN_i)\, \chrfcn{Y\ssN_i(t)\ge y}(\omega) - 
\int_{\msp} g(\rho)\,e^{-\rho((t-t_0(y,t),t])}\, \Lambda(d\rho) 
\biggr|
\\ \dsp {} \le 
\sup_{t\in[0,T];\ y_C(t)\ge y}
\biggl|\biggl(\frac1N\sum_{i=1}^N g(\rho\ssN_i)
-\int_{\msp} g(\rho)\, \Lambda(d\rho) \biggr)
\\ \dsp \phantom{{}\le}
-\biggl(\frac1N\sum_{i=1}^N g(\rho\ssN_i)\, \chrfcn{Y\ssN_i(t)< y}(\omega)
-\int_{\msp} g(\rho)\,(1-e^{-\rho((t-t_0(y,t),t])})\, \Lambda(d\rho) \biggr)
\biggr|
\\ \dsp {} \le \frac{5M}K\,,
\arre\end{equation}
for $\omega\in\tilde{\Omega}_K$ and $N>N_0(\omega)$.

Next, let $y_C(t)\le y$, and let $j=j(t)$ be the integer such that
\eqnb
\label{e:fnctHDLproof299}
y_{B,j}(t)\le y<y_{B,j+1}(t).
\end{equation}
With an argument similar as that below \eqnu{fnctHDLproof290},
such an integer $j=j(t)$ exists if $y_C(t)\le y$.
Since $y_B(y_0,t)$ is increasing in $y_0$\,, 
\eqnu{fnctHDLproof299} also implies 
\eqnb
\label{e:fnctHDLproof2990}
y_{0,j}<\hat{y}(y,t)\le y_{0,j+1}.
\end{equation}
Since \eqnu{fnctHDLproof24} implies
\[ 0\le y-y_{B,j}(t) \le y_{B,j+1}(t)-y_{B,j}(t)<\frac1K\,, \]
with \eqnu{fnctHDLproof28} and a similar argument as for \eqnu{HDLproof322}, 
we have
\eqnb
\label{e:fnctHDLproof293}
\arrb{l}\dsp
 \biggl| \frac1N\sum_{i=1}^N g(\rho\ssN_i)\, 
(\chrfcn{Y\ssN_i(t)\ge y} - \chrfcn{Y\ssN_i(t)\ge Y\ssN_B(y_{0,j},t)})
(\omega) \biggr|
\\ \dsp {}
\le \frac1N\sum_{i=1}^N |g(\rho\ssN_i)|\,  \biggl| 
\chrfcn{Y\ssN_i(t)\ge y} - \chrfcn{Y\ssN_i(t)\ge y_{B,j}(t)} \biggr|(\omega)
\\ \dsp \phantom{\le}
+  \frac1N\sum_{i=1}^N |g(\rho\ssN_i)|\, \biggl|
\chrfcn{Y\ssN_i(t)\ge y_{B,j}(t)} - \chrfcn{Y\ssN_i(t)\ge Y\ssN_B(y_{0,j},t)}
 \biggr|(\omega)
\\ \dsp {}
 \le M (y-y_{B,j}(t))+ M|Y\ssN_B(y_{0,j},t)(\omega)-y_{B,j}(t) |
<\frac{2M}K\,.
\arre
\end{equation}
Note also that, as in the argument for \eqnu{HDLproof321},
\eqnb
\label{e:fnctHDLproof294}
\chrfcn{(x\ssN_i-1)/N \ge y_{0,j},\ \tau\ssN_{i,1}>t}
= \chrfcn{Y\ssN_i(t)\ge Y\ssN_B(y_{0,j},t)}.
\end{equation}

Adding up \eqnu{fnctHDLproof293},
\eqnu{fnctHDLproof26} and \eqnu{fnctHDLproof29001},
and using  \eqnu{fnctHDLproof294} and triangular inequality, we arrive at
\eqnb
\label{e:fnctHDLproof11}
\arrb{l}\dsp
\sup_{t\in[0,T];\ y_C(t)\le y}\biggl|
\frac1N\sum_{i=1}^N g(\rho\ssN_i)\, \chrfcn{Y\ssN_i(t)\ge y}(\omega) -
\int_{\msp} g(\rho)\,e^{-\rho((0,t])}\,\mu_{0}(d\rho\times [\hat{y}(y,t),1))
\biggr|
\\ \dsp {}
\le \frac{4M}K\,,
\arre
\end{equation}
for $\omega\in\tilde{\Omega}_K$ and $N>N_0(\omega)$.

Combining \eqnu{fnctHDLproof1} and \eqnu{fnctHDLproof11}, we have
\[ \arrb{l}\dsp
\sup_{t\in[0,T]}
\biggl|\int_{\msp} g(\rho)\,\mu\ssN_t(d\rho\times [y,1))(\omega)
-\int_{\msp} g(\rho)\,\mu_t(d\rho\times [y,1))\biggr|\le\frac{5M}K\,,
\\ \dsp
\ N>N_0(\omega),\ \omega\in\tilde{\Omega}_K\,.
\arre\]
Finally, put $\tilde{\Omega}=\bigcap_{K=1}^{\infty}\tilde{\Omega}_K$.
Then $\prb{\tilde{\Omega}}=1$.
Let $\omega\in\tilde{\Omega}$.
For any $\eps>0$ take an integer $K$ such that $K>{5M}/{\eps}$\,.
Then $\omega\in \tilde{\Omega}\subset \tilde{\Omega}_K$ implies
\[
\sup_{t\in[0,T]}
\biggl|\int_{\msp} g(\rho)\,\mu\ssN_t(d\rho\times [y,1))(\omega)
-\int_{\msp} g(\rho)\,\mu_t(d\rho\times [y,1))\biggr|\le\frac{5M}K<\eps,
\]
for $N>N_0(\omega)$,
which implies \eqnu{fnctHDLproof0}, and therefore
\lemu{fnctBauer} implies the Theorem.
\prfe

\section{Case when the intensities have common time dependence.}
\seca{Appl}

To consider the case where the intensity measure $\rho$ has a density,
denote the set of locally integrable functions on $\preals$
by $\fsp$. $\fsp$ is a complete separable metric space.
Let $\iota$ be a map $\iota:\ \fsp\to\msp$ which maps $\tilde{w}\in\fsp$
to the measure on $\preals$ with density $\tilde{w}$ determined by
\eqnb
\label{e:iota}
 \iota(\tilde{w})((s,t])=\int_s^t \tilde{w}(u)\,du,\ 0\le s< t. 
\end{equation}
\prpb
\prpa{L1}
Assume that $\tilde{w}\ssN_i\in\fsp$, $i=1,2,\ldots,N$, $N=1,2,\ldots$, 
and for each $N$, put
\[ \tilde{\Lambda}\ssN=\frac1N\sum_{i=1}^N \delta_{\tilde{w}\ssN_i}\,. \]
If there exists a probability distribution $\tilde{\Lambda}$ on $\fsp$
such that $\tilde{\Lambda}\ssN$ converges weakly to $\tilde{\Lambda}$
as $N\to\infty$, then
the sequence of distribution $\Lambda\ssN$, $N=1,2,\ldots$,
on the set of intensity measures $\msp$ defined by
$\Lambda\ssN=\tilde{\Lambda}\ssN\circ \iota^{-1} $,
with $\iota$ as in \eqnu{iota},
converges weakly as $N\to\infty$ to
$\Lambda:=\tilde{\Lambda}\circ \iota^{-1}. $
Moreover, for all $0\le s<t$,
$\lambda\ssN_{s,t}$ defined by \eqnu{lambdaNstLambda}
converges weakly as $N\to\infty$ to $\lambda_{s,t}$ defined by
\eqnu{lambdaN2lambda}.
\DDD\prpe
\prfb
Let $g:\ \msp\to\reals$ be a bounded continuous function on $\msp$.
Then the definitions imply
\[
\int_{\msp} g(\rho)\,\Lambda\ssN(d\rho)
=\int_{\fsp} g(\iota(\tilde{w}))\,\tilde{\Lambda}\ssN(d\tilde{w}).
\]
Let $\{\tilde{w}_n\}$ be a sequence converging in $\fsp$ to $\tilde{w}$,
and let 
$f:\ \preals\to\reals$ be a continuous function with compact support:
$f(u)=0$, $u\ge k$, for some integer $k$.
Then $f$ is bounded: $|f(u)|\le M$, $u\in\preals$, for some $M$.
Hence
\[ \arrb{l}\dsp
\left|\int_{\preals} f(u)\, \tilde{w}_n(u)\,du
-\int_{\preals} f(u)\, \tilde{w}(u)\,du\right| 
=
\left|\int_0^k f(u)\, \tilde{w}_n(u)\,du
-\int_0^k f(u)\, \tilde{w}(u)\,du\right| 
\\ \dsp {}
\le M\,\int_0^k | \tilde{w}_n(u)-\tilde{w}(u)|\,du
\to 0, \ n\to\infty.
\arre \]
This holds for all continuous function $f$ with compact support,
hence $\lim_{n\to\infty} \iota(\tilde{w}_n)= \iota(\tilde{w})$ in vague 
topology, which further implies
\[\limf{n} g(\iota(\tilde{w}_n)) = g(\iota(\tilde{w})).\]
This proves that $g\circ \iota:\ \fsp\to \reals$ is a bounded continuous
function, hence the assumption $\tilde{\Lambda}\ssN\to\tilde{\Lambda}$
implies
\[
\limf{N} \int_{\msp} g(\rho)\,\Lambda\ssN(d\rho)
= \int_{\msp} g(\rho)\,\Lambda(d\rho).
\]
This holds for any bounded continuous function $g$,
which proves
$\Lambda\ssN\to \Lambda$, weakly as $N\to\infty$.

Let $t>s>0$ and put $b[\tilde{w}]=\int_s^t \tilde{w}(u)\,du$.
In a similar way as above, the definitions imply
\[
\lambda\ssN_{s,t}=
\int_{\fsp} \delta_{b[\tilde{w}]} \tilde{\Lambda}\ssN(d\tilde{w})
\ \mbox{ and }\ 
\lambda_{s,t}=
\int_{\fsp} \delta_{b[\tilde{w}]} \tilde{\Lambda}(d\tilde{w}). 
\]
Let $h:\ \preals\to\reals$ be a bounded continuous function.
Then the map
\[ \fsp\ni\tilde{w}\mapsto h(b[\tilde{w}])\in\reals\]
is bounded and continuous, hence the assumption
$\tilde{\Lambda}\ssN\to \tilde{\Lambda}$ implies
\[ \arrb{l}\dsp
\int_{\preals} h(w) \,\lambda\ssN_{s,t}(dw)=
\int_{\fsp} h(b[\tilde{w}]) \tilde{\Lambda}\ssN(d\tilde{w})
\\ \dsp {} \to 
\int_{\fsp} h(b[\tilde{w}])  \tilde{\Lambda}(d\tilde{w})
=\int_{\preals} h(w) \,\lambda_{s,t}(dw),\ N\to\infty,
\arre \]
hence 
$\lambda\ssN_{s,t}\to \lambda_{s,t}$, weakly as $N\to\infty$.
\prfe

\prpu{L1} implies that the assumption \eqnu{lambdaN2lambda} 
in \thmu{HDLinhomog} is redundant if the intensity measures have 
densities.

For the rest of this section, we further assume a common time dependence 
for all $\tilde{w}\ssN_i$ in \prpu{L1}. Namely, we assume that there exist
$\tilde{a}\in \fsp$ and positive constants 
\[ w\ssN_i>0,\ i=1,2,\ldots,N,\ N=1,2,\ldots, \]
such that the intensity measure of the Poisson random measures $\nu\ssN_i$
in the stochastic ranking process \eqnu{SIformSRP} is given by
\eqnb
\label{e:commontimedep}
\rho\ssN_i((s,t])=w\ssN_i\, \int_s^t \tilde{a}(u)\,du,
\ \ i=1,2,\ldots,N,\ N=1,2,\ldots.
\end{equation}
As in the proof of \prpu{L1}, we have 
\corb
\cora{R2L1}
Let $\tilde{a}\in \fsp$.
If there exists a probability distribution $\lambda$ on $\preals$
such that
\eqnb
\label{e:HDLlimitPareto}
\lambda\ssN:=\frac1N\sum_{i=1}^N \delta_{w\ssN_i}\to\lambda,
\ \mbox{ weakly, as }\ N\to\infty,
\end{equation}
then a sequence of probability distributions 
$\tilde{\Lambda}\ssN$, $N=1,2,\ldots$, on $\fsp$ defined by
\[
\tilde{\Lambda}\ssN=\int_{\preals} \delta_{w\tilde{a}}\lambda\ssN_i(dw)
=\frac1N\sum_{i=1}^N \delta_{w\ssN_i\tilde{a}}
\]
converges weakly to a probability distribution
$\tilde{\Lambda}=\int_{\preals} \delta_{w\tilde{a}} \lambda(dw)$,
as $N\to\infty$.

In particular, \prpu{yC} holds with
$\rho\ssN_i((s,t])=w\ssN_i \int_s^t \tilde{a}(u)\,du$,
and $y_C(t)$ of \eqnu{yC} is given by
\eqnb
\label{e:yCx}
y_C(t)=1-\int_{\preals} e^{-w\, A(t)}\,\lambda(dw),
\end{equation}
where
\eqnb
\label{e:A}
A(t)=\int_0^t \tilde{a}(u)\,du.
\end{equation}
\DDD\core

The formula \eqnu{yCx} is to be compared with the case
of the (homogeneous) Poisson process in \cite[Proposition 2]{HH071},
where we have
\eqnb
\label{e:yCt}
y_C(t)=1-\int_{\preals} e^{-wt} \lambda(dw).
\end{equation}
$\lambda$ in \eqnu{yCt}
is the (infinite particle limit asymptotic) distribution of 
jump rates, while $\lambda$ in the case of common time dependence
\eqnu{yCx} is the distribution of relative jump rates.

To study a time change according to the common intensity measure,
let us first make a heuristic observation.
Suppose we could trace the trajectories of  $n\le N$ particles 
$j_1,j_2,\ldots,j_n$.
The total number of jumps of the $n$ particles in the time interval $(0,t]$
is given by
\eqnb
\label{e:ntotaljump}
S\ssNn(t)=\sum_{i=1}^n \nu\ssN_{j_i}((0,t]).
\end{equation}
If $n$ is large ($n\gg1$), we expect
as a consequence of the law of large numbers,
as in \prpu{yC},
\eqnb
\label{e:SANn}
 S\ssNn(t)\simeq \sum_{i=1}^n \rho\ssN_{j_i}((0,t])=A(t)\,Z(N,n),
\end{equation}
where we put
\eqnb
\label{e:totaljumpsscalingNn}
Z(N,n)=\sum_{i=1}^n w\ssN_{j_i}\,,
\end{equation}
and also used \eqnu{commontimedep} and \eqnu{A}.
Using \eqnu{SANn} in \eqnu{yCx},
we have
\eqnb
\label{e:timechangeyCNn}
y_C(t)\simeq 1-\int_{\preals} e^{\dsp -w \,  S\ssNn(t) /Z(N,n)}\,\lambda(dw).
\end{equation}
The approximate formula \eqnu{timechangeyCNn} suggests that, 
if we perform a time change $t'=S\ssNn(t)$, 
then modulo scaling constant $Z(N,n)$,
we recover a formula \eqnu{yCt} for the homogeneous case.

We can put the heuristic consideration which lead to
\eqnu{timechangeyCNn} in a mathematically precise form. 
For $t\ge0$, let
\eqnb
\label{e:totaljump}
 S\ssN(t)=\sum_{i=1}^N \nu\ssN_i((0,t])
\end{equation}
and denote its right continuous inverse by
\eqnb
\label{e:totaljumptimechange}
s\ssN(t)=\inf\{s\ge0\;;\ S\ssN(s)>t \}.
\end{equation}
Let $\tilde{a}\in \fsp$.
For simplicity, assume further that
\eqnb
\label{e:posdensityintensity}
\tilde{a}(t)>0,\ t\ge0.
\end{equation}
Then $A(t)$ of \eqnu{A} is strictly increasing,
and the inverse function $A^{-1}$ is also continuous.
\thmb
\thma{yCttimechange}
Let $\tilde{a}\in \fsp$, and assume \eqnu{posdensityintensity}.
Put
\eqnb
\label{e:totaljumpsscaling}
Z(N)=\sum_{i=1}^N w\ssN_i
\end{equation}
and assume
\eqnb
\label{e:totaljumpsscalinginfinity}
\limf{N} Z(N)= \infty.
\end{equation}
If, as in \coru{R2L1},
there exists a probability distribution $\lambda$ on $\preals$
such that \eqnu{HDLlimitPareto} holds,
then for each $t\ge0$
\eqnb
\label{e:yCtcommontimedepscaledbytotaljumps}
Y\ssN_C(s\ssN(Z(N)\,t)) \to y_C(A^{-1}(t))
= 1-\int_{\preals} e^{-w\,t} \lambda(dw),
\ \mbox{ in probability, as }\ N\to\infty,
\end{equation}
where $Y\ssN_C$ is defined in \eqnu{YCN}.
\DDD\thme
To prove \thmu{yCttimechange}, we first provide a rigorous version of
\eqnu{SANn}.
\lemb
\lema{commontimedepscaledtotaljumpsLLN}
For $t\ge0$,
\eqnb
\label{e:commontimedeptotaljumpsLLN}
 \frac1{Z(N)} S\ssN(t)\to A(t),\ \mbox{ in probability, as }\ N\to\infty.
\end{equation}
and
\eqnb
\label{e:commontimedepscaledtotaljumpsLLN}
 s\ssN(Z(N)\,t)\to A^{-1}(t),\ \mbox{ in probability, as }\ N\to\infty.
\end{equation}
\DDD\leme
\prfb
Since by definition $\nu\ssN_i((0,t])$ follows the Poisson distribution
with expectation $\rho\ssN_i((0,t])$, we have
\eqnb
\label{e:commontimedeptotaljumpsmean}
\EE{S\ssN(t)}=\VV{S\ssN(t)}=A(t)\, Z(N),
\end{equation}
where $\VV{\cdot}$ denotes variance.
For $\eps>0$, \eqnu{commontimedeptotaljumpsmean}, \eqnu{totaljumpsscaling},
and Chebyshev's inequality imply
\[
\prb{|S\ssN(t)-\EE{S\ssN(t)}|>Z(N)\eps}\le (\eps Z(N))^{-2}\VV{S\ssN(t)}
= \frac{A(t)}{\eps^2 Z(N)}\,,
\]
which, with \eqnu{totaljumpsscalinginfinity}, implies
\[
\frac1{Z(N)} (S\ssN(t)-\EE{S\ssN(t)})\to 0,\ \mbox{ in probability, as }
\ N\to\infty.
\]
This, with \eqnu{commontimedeptotaljumpsmean}, 
implies \eqnu{commontimedeptotaljumpsLLN}.

Next, noting that $S\ssN(t)$ is non-decreasing in $t$,
\eqnu{totaljumptimechange} implies
\eqnb
\label{e:commontimedepscaledtotaljumpsLLN1}
\{\omega\in\Omega\;;\ s\ssN(Z(N)t)(\omega) \ge A^{-1}(t)+\eps\}
\subset \{\omega\in\Omega\;;\ 
 \frac1{Z(N)} S\ssN\left(A^{-1}(t)+\frac{\eps}2\right)(\omega)\le t\}.
\end{equation}
The assumption \eqnu{posdensityintensity} implies
that $A$ is strictly increasing,
hence, 
$\delta= A(A^{-1}(t)+{\eps}/2)-t > 0$,
and
\[\arrb{l}\dsp
\{\omega\in\Omega\;;\ \frac1{Z(N)} S\ssN(A^{-1}(t)+\frac{\eps}2)\le t \}
\\ \dsp {}
=
\{\omega\in\Omega\;;\ \frac1{Z(N)} S\ssN(A^{-1}(t)+\frac{\eps}2)\le 
A\left(A^{-1}(t)+\frac{\eps}2\right)-\delta \}
\\ \dsp {}
\subset
\{\omega\in\Omega\;;\ 
\left| \frac1{Z(N)} S\ssN(A^{-1}(t)+\frac{\eps}2)
-A\left(A^{-1}(t)+\frac{\eps}2\right)\right|
\ge \delta \}.
\arre\]
This and \eqnu{commontimedeptotaljumpsLLN} and
\eqnu{commontimedepscaledtotaljumpsLLN1} imply
\eqnb
\label{e:commontimedepscaledtotaljumpsLLN3}
\limf{N} \prb{s\ssN(Z(N)t) \ge A^{-1}(t)+\eps}=0.
\end{equation}
Similarly,
$
\delta' = t-A(A^{-1}(t)-{\eps}/2)>0$, and
\[\arrb{l}\dsp
\{\omega\in\Omega\;;\  s\ssN(Z(N)t) \le A^{-1}(t) -\eps \}
\\ \dsp {}
\subset \{\omega\in\Omega
\;;\ \frac1{Z(N)} S\ssN\left(A^{-1}(t)-\frac{\eps}2\right)\ge t \}
\\ \dsp {}
\subset
\{\omega\in\Omega\;;\ 
\left| \frac1{Z(N)} S\ssN(A^{-1}(t)-\frac{\eps}2)
-A\left(A^{-1}(t)-\frac{\eps}2\right)\right|
\ge \delta' \},
\arre \]
which implies
\eqnb
\label{e:commontimedepscaledtotaljumpsLLN4}
\limf{N} \prb{s\ssN(Z(N)t) \le A^{-1}(t)-\eps}=0.
\end{equation}
\eqnu{commontimedepscaledtotaljumpsLLN3} and
\eqnu{commontimedepscaledtotaljumpsLLN4} prove
\eqnu{commontimedepscaledtotaljumpsLLN}.
\prfe

\prfofb{\protect\thmu{yCttimechange}}
By triangular inequality, we have
\[ \arrb{l}\dsp
| Y\ssN_C(s\ssN(Z(N)\,t)) - y_C(A^{-1}(t)) |
\\ \dsp {}
\le 
| Y\ssN_C(s\ssN(Z(N)\,t)) -  Y\ssN_C(A^{-1}(t))|
+ |Y\ssN_C(A^{-1}(t)) - y_C(A^{-1}(t)) |.
\arre \]
\coru{R2L1} implies that
the second term in the right hand side converges to $0$ in probability
as $N\to\infty$, so it suffices to prove that, for all $\eps>0$,
\eqnb
\label{e:yCttimechangeprf2}
\limf{N} \prb{| Y\ssN_C(s\ssN(Z(N)\,t)) -  Y\ssN_C(A^{-1}(t))|\ge \eps}=0
\end{equation}
holds.

For $\delta>0$ put
\eqnb
\label{e:yCttimechangeprf3}
 \Omega\ssN_{\delta}:=\{\omega\in\Omega\;;\ 
|s\ssN(Z(N)\,t)(\omega) - A^{-1}(t)| < \delta\}.
\end{equation}
Then \eqnu{commontimedepscaledtotaljumpsLLN} implies
\eqnb
\label{e:yCttimechangeprf1}
\limf{N} \prb{ \Omega\ssN_{\delta}\,{}^c}=0.
\end{equation}
The definition \eqnu{YCN} of $Y\ssN_C$ implies 
\eqnb
\label{e:yCttimechangeprf4}
\arrb{l}\dsp
| Y\ssN_C(s\ssN(Z(N)\,t)) -  Y\ssN_C(A^{-1}(t))|
\\ \dsp {}
=\frac1N \sum_{i=1}^N \chrfcn{ s\ssN(Z(N)\,t)<\tau\ssN_{i,1} \le A^{-1}(t) }
+\frac1N \sum_{i=1}^N \chrfcn{ A^{-1}(t)<\tau\ssN_{i,1} \le s\ssN(Z(N)\,t) }.
\arre
\end{equation}
Combining \eqnu{yCttimechangeprf3} and \eqnu{yCttimechangeprf4}, we have
\[
 \prb{ | Y\ssN_C(s\ssN(Z(N)\,t)) -  Y\ssN_C(A^{-1}(t))| \ge \eps,\
\Omega_{\delta} }
\le
 \mathrm{P}\left[\;
  \sum_{i=1}^N
 \chrfcn{ \tau\ssN_{i,1} \in (A^{-1}(t)-\delta, A^{-1}(t)+\delta) }
 \ge N\eps
\;\right].
\]
Applying Chebyshev's inequality, we further have
\[\arrb{l}\dsp
 \prb{ | Y\ssN_C(s\ssN(Z(N)\,t)) -  Y\ssN_C(A^{-1}(t))| \ge \eps,\
\Omega_{\delta} }
\\ \dsp {}
\le
\frac1{N\eps} \sum_{i=1}^N 
 \EE{\chrfcn{ \tau\ssN_{i,1} \in (A^{-1}(t)-\delta, A^{-1}(t)+\delta) }}
\\ \dsp {}
=
\frac1{N\eps} \sum_{i=1}^N 
\left(
 e^{ -A(A^{-1}(t)-\delta)\, w\ssN_i }- e^{ -A(A^{-1}(t)+\delta)\, w\ssN_i }
\right)
\\ \dsp {}
=
\frac1{\eps}\int_{\preals}
\left(
 e^{ -A(A^{-1}(t)-\delta)\, w }- e^{ -A(A^{-1}(t)+\delta)\, w }
\right) \, \lambda\ssN(dw).
\arre \]
This, with \eqnu{yCttimechangeprf1} and 
the assumption \eqnu{HDLlimitPareto}, implies
\[ \arrb{l}\dsp
\limsup_{N\to\infty}
 \prb{| Y\ssN_C(s\ssN(Z(N)\,t)) -  Y\ssN_C(A^{-1}(t))|\ge \eps}
\\ \dsp {}
\le
\limsup_{N\to\infty} \prb{ {\Omega\ssN_{\delta}}^c}
+
\limsup_{N\to\infty}
 \prb{| Y\ssN_C(s\ssN(Z(N)\,t)) -  Y\ssN_C(A^{-1}(t))|\ge \eps,
\ \Omega\ssN_{\delta}}
\\ \dsp {}
\le
\frac1{\eps}\int_{\preals}
\left(
 e^{ -A(A^{-1}(t)-\delta)\, w }- e^{ -A(A^{-1}(t)+\delta)\, w }
\right) \, \lambda(dw).
\arre \]
This holds for all $\delta>0$, hence the bounded convergence theorem
and the continuity of $A(t)$ imply
\[ \arrb{l}\dsp
\limsup_{N\to\infty}
 \prb{| Y\ssN_C(s\ssN(Z(N)\,t)) -  Y\ssN_C(A^{-1}(t))|\ge \eps}
\\ \dsp {}
\le 
\inf_{\delta>0}
\frac1{\eps}\int_{\preals}
\left(
 e^{ -A(A^{-1}(t)-\delta)\, w }- e^{ -A(A^{-1}(t)+\delta)\, w }
\right) \, \lambda(dw)
\\ \dsp {}
\le 
\frac1{\eps}\int_{\preals}
\lim_{\delta\downarrow0}
\left(
 e^{ -A(A^{-1}(t)-\delta)\, w }- e^{ -A(A^{-1}(t)+\delta)\, w }
\right) \, \lambda(dw)
=0.
\arre \]
This proves \eqnu{yCttimechangeprf2},
hence \thmu{yCttimechange} is proved.
\QED\prfofe

As an explicit example to $Z(N)$ and $\lambda$, consider, as in 
\cite{HH072,HH073}, the Zipf's law, which is
\eqnb
\label{e:Zipf}
w\ssN_i = a \left(\frac{N}{i}\right)^{1/b},\ \ i=1,2,\ldots,N,
\end{equation}
for positive constants $a$ and $b$.
For this choice,
\eqnb
\label{e:totaljumpsscalingZipf}
Z(N)=\sum_{i=1}^N w\ssN_i
 = (1+o(N^0))\times \left\{ \arrb{ll}
\dsp aN\int_0^1 x^{-1/b} dx=\frac{aN\,b}{b-1} & b>1, \\ 
\dsp aN\int_{1/N}^1 x^{-1} dx= aN \log N & b=1, \\
\dsp aN^{1/b}\sum_{i=1}^{\infty} \frac1{i^{1/b}}
= aN^{1/b}\zeta(1/b) & 0<b<1\,.
\arre\right.
\end{equation}
The corresponding $N\to\infty$ weak limit is
the (generalized) Pareto distribution, defined by
\eqnb
\label{e:Paretolambda}
\lambda([w,\infty))= \left\{ \arrb{ll} \dsp \left(\frac{a}{w}\right)^b
 & w\ge a, \\ 1 & w<a. \arre \right.
\end{equation}
With the Pareto distribution \eqnu{Paretolambda} for $\lambda$, 
\eqnu{timechangeyCNn} is (for $N=n$)
\eqnb
\label{e:timechangexCNPareto}
\arrb{l}\dsp
x_C(t)=Ny_C(t)+1\simeq
 N-N\int_{\preals} e^{\dsp -w \,  S\ssN(t) /Z(N)}\,\lambda(dw)
\\ \dsp {}
=N-b\,\left(\frac{S\ssN(t)}{\zeta_N(1/b)}\right)^b\,
 \Gamma(-b,\frac{S\ssN(t)}{N^{1/b}\zeta_N(1/b)})
\\ \dsp {}
=N-Ne^{-S\ssN(t)/(N^{1/b}\zeta_N(1/b))}
+\left(\frac{S\ssN(t)}{\zeta_N(1/b)}\right)^b\, 
 \Gamma(1-b,\frac{S\ssN(t)}{N^{1/b}\zeta_N(1/b)}) =: x\ssN_b(S\ssN(t)),
\arre
\end{equation}
where
$\zeta_N(z)= \sum_{i=1}^{N} i^{-z}$.
The last line in \eqnu{timechangexCNPareto} is obtained by 
integration by parts from the second line, as in \cite{HH073},
and is suitable for $0<b<1$.
Note that the parameter $a$ in the Pareto distribution \eqnu{Paretolambda}
disappears in the time changed formula \eqnu{timechangexCNPareto}.

\appendix

\section{Remarks on practical application.}

In \cite{HH072,HH073}, the mathematical results on the stochastic
ranking processes has been successfully applied to practical data,
such as ranking data of books at an online bookstore Amazon.co.jp
\cite{HH073,HH072} and list of subject titles at a collected 
bulletin board 2ch.net \cite{HH072}.

One may wonder why such a simple rule as the move-to-front rule
could be observed in actual social activities. An explanation is that the
ranking numbers on the web (such as those representing the books, 
in the case of online bookstores) usually seek to align the web pages
in the order of \textit{current popularity} of the pages.
A social impact of the development of web-based activities
is that it has become possible to catalog a huge amount of unpopular items
\cite{longtail}.
In fact, a majority of books catalogued on an online bookstore are
sold less than one copy a month. For such books, any reasonable 
order reflecting the current popularity would be equal to the order
of the time of most recent sales, because the second recent sale of such
book would be long ago, hence would not reflect current popularity.
Thus the move-to-front rule will provide a simple but \textit{universal} model
in the rankings on the web.

A ranking of a book at Amazon.co.jp jumps close to top of the ranking
whenever the book is sold at Amazon.co.jp \cite{HH073},
and a subject title in the web page for the list of 2ch.net 
jumps to the top whenever a comment (a `response') concerning the subject
is written \cite{HH072}. Ordering a book and responding to a subject
are social activities which naturally are expected to contain 
day-night difference in the intensity.

Explicit time dependence,
reflecting day-night difference of social activities, 
are observed in actual data.
Let us  regard such time dependence as the non-uniformity of 
intensity measures $\rho\ssN_i$.
$\rho\ssN_i$ are usually unknown quantities to be determined statistically 
from observed data.
We then have to consider both particle dependence and
time dependence in the statistical analysis of the practical data.
The assumption of common time dependence \eqnu{commontimedep} developed in
\secu{Appl} provides a simple way to 
take day-night-difference of social activity into account,
in applying the stochastic ranking process with inhomogeneous intensity.

\subsection{Factorization of day-night social activity difference.}

In \cite{HH073}, 
a data taken during the period of about 3 months at Amazon.co.jp
is used to statistically obtain $\lambda$, based on \eqnu{yCt}.
The data was taken manually in the year 2007, at 21:00 each day.
We can show that in the case of common time dependence assumption
\eqnu{commontimedep}, we can `factorize' periodic time dependence of 
$\tilde{a}$, and that the use of \eqnu{yCt} in \cite{HH073,HH072} is 
justified in obtaining $\lambda$ from data with periodic time dependence.
In fact, assume that there exists a positive constant $T$ such that
\eqnb
\label{e:periodictildea}
\tilde{a}(t+T)=\tilde{a}(t),\ t\ge0.
\end{equation}
We may normalize $w\ssN_i$'s in \eqnu{commontimedep} so that
\eqnb
\label{e:normalizedtildea}
\frac1T\int_0^{T} \tilde{a}(u)\,du=1
\end{equation}
holds. Then \eqnu{periodictildea} and \eqnu{normalizedtildea} imply
$\int_{t}^{t+T} (\tilde{a}(u)-1)\,du =0$, so that
\eqnb
\label{e:ap}
A_p(t):= A(t)-t =\int_0^t (\tilde{a}(u)-1)\,du
\end{equation}
is a periodic function with period $T$, and \eqnu{yCx} is
\eqnb
\label{e:yCxAp}
y_C(t)=1-\int_{\preals} e^{-w\, (t+A_p(t))}\,\lambda(dw).
\end{equation}
If we collect data at each fixed time of the day, at
$t_n=t_0+n\,T$, $n=0,1,2,\ldots$, then \eqnu{yCxAp} implies
\eqnb
\label{e:yCtn}
y_C(t_n)=1-\int_{\preals} e^{-w (nT+ t_0+A_p(t_0))}\,\lambda(dw).
\end{equation}
Hence the effect of day-night difference in $\tilde{a}$ is absorbed
in the translation of origin of time $t_0\mapsto t_0+A_p(t_0)$,
and the use of formula \eqnu{yCt} for the constant intensity is justified.

A consideration of this subsection is of practical use 
when one has a data much longer than $24$ hours,
as in the case of \cite{HH073}.

\subsection{Time change according to intensity measure.}

In \cite{HH072}, a data of list of subject (`thread') titles at a collected 
bulletin board 2ch.net is statistically analyzed using 
stochastic ranking process.
In \cite{HH072} the data was collected from a short period
in the daytime, and the problem
of day-night activity difference was not serious,
hence a fit to the formula \eqnu{yCt} for the constant jump rate
(homogeneous intensity) was possible \cite{HH072}.
However, to study data of longer periods for sharper statistical results,
effects of day-night activity difference need to be taken into account.

In applying \eqnu{commontimedep} to the obtained data
to extract time dependence (day-night difference),
we need to estimate the function $\tilde{a}$ in \eqnu{commontimedep}
or $A$ in \eqnu{A}.
This is accomplished by making use of \eqnu{SANn} and \eqnu{timechangeyCNn}.
In the case of 2ch.net \cite{HH072}, 
$N$ in \eqnu{SANn} or \eqnu{timechangeyCNn} is about $700$, 
and since full records of transaction are accessible at 2ch.net,
it is possible to put $n$ in \eqnu{timechangeyCNn} equal to $N$
and count all the threads' jumps.
In the case of Amazon.co.jp,
$N$ is of order million, and $n=N$ is unrealistic.
Even in such cases, if we observe sufficiently large number of books
($n\gg 1$), we can apply the idea introduced here.

Note that the series $Z(N)$ are approaching 
their asymptotics in \eqnu{totaljumpsscalingZipf} rather slowly
for the Pareto distribution.
Therefore in practical application of 
\eqnu{yCtcommontimedepscaledbytotaljumps}
with the Pareto distribution for $\lambda$,
if one takes $N=O(10^3)$ as in 2ch.net \cite{HH072},
one should avoid using the asymptotic formula in
the right hand side of \eqnu{totaljumpsscalingZipf},
and calculate the finite sums \eqnu{totaljumpsscaling} or 
\eqnu{totaljumpsscalingNn}.

We announce that we actually collected a $24$ hours data of size $n_d=70140$
from 2ch.net, 
and performed a statistical fit of the data to \eqnu{timechangexCNPareto},
with $N=697$, and obtained $b=0.872\pm0.002$.
(The error is 90\% confidence level. See \cite{HH073} for details.)
Apparently, we have a good single parameter fit to the data,
which suggests that the practical assumption \eqnu{commontimedep}
is good. 
Details may be reported elsewhere.

We note that in \cite{HH072}, a value of $b=0.6145$ was obtained for
2ch.net (with different set of data). This is much smaller than the
present result. The data used in \cite{HH072} was small in size,
because the data was collected manually in those times,
and also, to avoid influence of day-night difference in the total activity,
the data was for a short time period in \cite{HH072}, so that
the result in \cite{HH072} is less reliable compared to the present result.

We also note that we have $b<1$, consistently with previous observation
\cite{HH073} for Amazon.co.jp, where we obtained $b=0.809$.
This shows that, as in Amazon.co.jp, 
the popularity of subjects is concentrated to a relatively small 
number of threads in 2ch.net.

\address{
Yuu Hariya \\
Mathematical Institute \\
Tohoku University \\
Sendai 980-8578 \\
Japan
}{hariya@math.tohoku.ac.jp}
\address{
Kumiko Hattori \\ 
Department of Mathematics and Information Sciences \\
Tokyo Metropolitan University \\
Hachioji 192-0397 \\
Japan
}{khattori@tmu.ac.jp}
\address{
Tetsuya Hattori \\ 
Faculty of Economics \\
Keio University \\
Yokohama 223-8521 \\
Japan
}{
hattori@econ.keio.ac.jp}
\address{
Yukio Nagahata \\
Department of Mathematical Science \\
Osaka University \\
Toyonaka 560-8531 \\
Japan
}{nagahata@sigmath.es.osaka-u.ac.jp}
\address{
Yuusuke Takeshima \\ 
Fukoku Mutual Life Insurance Company \\
Tokyo 100-0011 \\
Japan
}{\scriptsize ds-y.825533\_p@agate.plala.or.jp}
\address{
Takahisa Kobayashi \\ 
Miyagiken Sendai Daini High School \\
Sendai 980-8631 \\
Japan
}{\scriptsize bassuy-kobakoba0911@yahoo.co.jp}

\enddocument